\documentclass[11pt]{amsart}
\usepackage{fullpage}
\usepackage{amssymb}
\usepackage{amstext}
\usepackage{amsmath}
\usepackage{amsthm}
\usepackage{stmaryrd}
\usepackage{verbatim}
\usepackage{comment}
\usepackage{hyperref}
\usepackage{tikz}

\newtheorem{theorem}{Theorem}[section]
\newtheorem{corollary}[theorem]{Corollary}
\newtheorem{proposition}[theorem]{Proposition}
\newtheorem{lemma}[theorem]{Lemma}
\theoremstyle{definition}
\newtheorem{definition}[theorem]{Definition}
\theoremstyle{remark}
\newtheorem{remark}[theorem]{Remark}
\newtheorem{example}[theorem]{Example}

\newcommand{\C}{\mathbb C}
\newcommand{\R}{\mathbb R}
\newcommand{\Q}{\mathbb Q}
\newcommand{\Z}{\mathbb Z}
\newcommand{\HH}{\mathfrak H}
\newcommand{\SL}{\mathrm{SL}}
\newcommand{\ind}{\mathbf 1}
\DeclareMathOperator{\ord}{ord}
\DeclareMathOperator{\rank}{rank}
\DeclareMathOperator{\Span}{span}

\title{Modular \texorpdfstring{$L$}{L}-functions and Chebyshev systems}
\author{Wei Wang}
\address{Department of Mathematics, Shaoxing University, Shaoxing 312000, China}
\email{weiwang\_math@163.com}

\subjclass[2020]{Primary 11F67; Secondary 11F11, 11G18, 41A50}
\keywords{Modular forms, completed $L$-functions, Chebyshev spaces, periods,
Atkin--Lehner involutions, real zeros}

\begin{document}

\begin{abstract}
We study real zeros of completed modular $L$-functions and the determination of cusp forms by their periods. Using Chebyshev systems and the Riemann–Roch theorem, we obtain bounds for real zeros, counted with multiplicity, and rank estimates for evaluation maps in terms of the dimensions of cusp-form spaces and the genera of modular curves. These estimates give criteria for finite families of periods to determine a cusp form and yield explicit dual bases in cases where the relevant quotient curves have genus zero. We also bound the dimensions of rational spans of eigenform periods.
\end{abstract}
\maketitle

\section{Introduction and main results}\label{section_intro}
Let $N\ge1$, let $k\ge2$ be even, and let $f(z)=\sum_{n\ge1}a_f(n)e^{2\pi inz}$ belong to $S_k(\Gamma_0(N))$. Define
\[
L(f,s)=\sum_{n\ge1}\frac{a_f(n)}{n^s},\qquad \Lambda_N(f,s)=N^{s/2}(2\pi)^{-s}\Gamma(s)L(f,s).
\]
The Dirichlet series initially converges in a right half-plane, and $\Lambda_N(f,s)$ extends to an entire function. Its critical values are, up to explicit nonzero factors, the periods
\begin{equation}\label{equation_period-definition}
r_n(f)=\int_0^{i\infty}f(z)z^n\,dz=i^{n+1}N^{-(n+1)/2}\Lambda_N(f,n+1),
\qquad 0\le n\le k-2.
\end{equation}
At level one, these periods satisfy the Eichler--Shimura relations. A fundamental question in period theory is to identify which periods suffice to determine a cusp form. The even and odd parts of the period polynomial each determine the cusp form \cite{kohnen1984modular}. Using the Eichler--Shimura relations, Fukuhara \cite[Theorem~2.2]{fukuhara2007explicit} constructed an explicit basis of the dual space consisting of odd periods, with further applications to Hecke operators. Ni et al.~\cite{ni2026explicit} used the same relations to construct explicit spanning families of periods and derive spanning results for Rankin--Cohen brackets.

Lei, Ni, and Xue \cite{lei2023lineara,lei2023linearb} studied period independence through products and Rankin--Cohen brackets of Eisenstein series and estimates for their Fourier coefficients. Their results cover up to three periods of the same parity, subject to dimension and index conditions, and any fixed number of such periods in sufficiently large weight. Xue \cite{xue2024linear} treated pairs of opposite parity under additional index restrictions. Refining Fourier coefficient estimates for the corresponding kernels, Ni and Xue \cite{ni2026periods} obtained the bound $k\gg n^2$ for certain initial families of $n$ odd or even periods. The present paper treats period independence through real zeros of Mellin transforms, allowing both the evaluation parameters and the level to vary.

Chebyshev systems, the main analytic tool in this paper, are classical in approximation theory, interpolation, and moment problems. They generalize the zero-counting and interpolation properties of algebraic polynomials \cite{borwein1995polynomials,karlin1966tchebycheff}. They underlie uniqueness results for best uniform approximation, generalized Gaussian quadrature, and the geometry of moment spaces \cite{karlin1966tchebycheff}. They also provide zero bounds for Abelian integrals, with applications to limit cycles bifurcating under perturbations of planar Hamiltonian systems \cite{grau2011chebyshev}.

The Mellin representation of $\Lambda_N(f,s)$ connects period independence with the zero-counting properties of Chebyshev systems. The Riemann--Roch theorem bounds the zeros of cusp forms along suitable paths on modular curves, and the sign-change principle transfers this control to their Mellin transforms. The resulting estimates apply to arbitrary real evaluation points in the stated intervals. Counting zeros with multiplicity also allows derivative evaluations to be treated by the same argument. These bounds hold in every admissible weight, with an optimal leading dimension term at fixed level, while the genus of the relevant modular curve controls the possible loss of interpolation rank. This gives a common geometric and analytic approach to real zeros, period independence, and interpolation.

For a nonzero analytic function $G$, write $Z_I(G)$ for its number of real zeros in $I$, counted with analytic multiplicity, including an endpoint if it belongs to $I$ and $G$ is analytic near it. Write $[x]_+=\max\{x,0\}$ and let $\ind_{\{\mathcal P\}}$ be $1$ when $\mathcal P$ holds and $0$ otherwise.

For a group $\Gamma$ and a character $\chi$, let $g(\Gamma)$ be the genus of the compactified quotient curve $X_\Gamma$, and let
\[
d_k(\Gamma,\chi)=\dim_{\C}S_k(\Gamma,\chi).
\]
We omit $\chi$ when it is trivial. For $N\ge1$, write
\[
g_N=g(X_0(N)),\qquad d_k(N)=\dim_{\C}S_k(\Gamma_0(N)).
\]
Dimensions of cusp-form spaces and ranks of their period maps are over $\C$ unless specified otherwise. We first state the results for $\Gamma_0(N)$.

\begin{theorem}\label{theorem_basic}
Let $N\ge1$ and let $k\ge2$ be even, with $d_k(N)>0$. Then every nonzero $f\in S_k(\Gamma_0(N))$ satisfies
\[
Z_{(-\infty,k/2)}(\Lambda_N(f,\cdot))\le d_k(N)+g_N-1-\ind_{\{k=2\}}.
\]
For any distinct real numbers $s_1,\ldots,s_q<k/2$,
\[
\begin{aligned}
&\rank\bigl(f\longmapsto(\Lambda_N(f,s_1),\ldots,\Lambda_N(f,s_q))\bigr)\\
&\hspace{12mm}\ge\min\{d_k(N),[q-g_N+\ind_{\{k=2\}}]_+\}.
\end{aligned}
\]
Thus a sufficient condition for injectivity of the evaluation map is $q\ge d_k(N)+g_N-\ind_{\{k=2\}}$. For $N>1$, both assertions hold on $\mathbb R$.
\end{theorem}
At level one, the evaluations at any $d_k(1)$ distinct real points below $k/2$ therefore form a dual basis. In particular, any set of critical periods $r_n$ with $n+1<k/2$ is linearly independent whenever its cardinality does not exceed $d_k(1)$, with no restriction on the parity of the selected indices.

The rank estimate can also be expressed in terms of matrices of completed $L$-values. Let $f_1,\ldots,f_{d_k(N)}$ be a basis of $S_k(\Gamma_0(N))$. For distinct real points $s_1,\ldots,s_q$, with $s_i<k/2$ when $N=1$, a sufficient condition for the evaluation matrix $\bigl(\Lambda_N(f_j,s_i)\bigr)_{1\le i\le q,\;1\le j\le d_k(N)}$
to have full column rank is $q\ge d_k(N)+g_N-\ind_{\{k=2\}}$. When $g_N=\ind_{\{k=2\}}$, every such square matrix is nonsingular.

The dimension term in Theorem \ref{theorem_basic} is necessary. Write $S_k(\Gamma_0(N))_{\R}$ for the space of cusp forms with real Fourier coefficients. Any $d_k(N)-1$ homogeneous real evaluation conditions have a nonzero solution in this space. Consequently, for $N>1$
and even $k\ge4$ with nonzero cusp-form space,
\begin{equation}\label{equation_dimension-sharpness}
\begin{aligned}
d_k(N)-1
&\le \max_{\substack{0\ne f\in S_k(\Gamma_0(N))_{\R}}}Z_{\R}(\Lambda_N(f,\cdot))\\
&\le d_k(N)+g_N-1.
\end{aligned}
\end{equation}
The zeros in the lower-bound argument may be prescribed at any $d_k(N)-1$ distinct real points. At fixed level, the gap in \eqref{equation_dimension-sharpness} is at most $g_N$, independently of the weight, so the leading dimension term is optimal. Corollary~\ref{corollary_rational-sharpness} further shows that in the genus-zero case, the upper bound can be attained by forms with rational Fourier coefficients, with simple zeros in prescribed disjoint intervals. The upper bound is also attained in positive-genus cases. In particular, Theorem~\ref{theorem_fricke-sharpness} shows that this holds whenever $X_0(N)$ has genus one. For each fixed even weight $k\ge2$, Corollary~\ref{corollary_et-finiteness} asserts that only finitely many levels $N>1$ with nonzero cusp-form space have the following property: the Mellin transform of every nonzero $f\in S_k(\Gamma_0(N))_{\R}$ has at most $d_k(N)-1$ real zeros, counted with multiplicity. The obstruction comes from the fixed points of the Fricke involution on $X_0(N)$: sufficiently many such points force some Mellin transforms to exceed this zero bound.
 
Our method also applies to other modular groups. Theorem~\ref{theorem_modular} extends these estimates to $\Gamma_1(N)$, $\Gamma(N)$, and full simultaneous Atkin--Lehner eigenspaces, with the genus of the corresponding quotient curve. Rank bounds for evaluations of $\Lambda_N(f,s)$ and its derivatives appear in Theorem~\ref{theorem_derivatives}, including the central derivatives allowed by the functional equation. For spaces with a fixed sign in the functional equation, these estimates also bound the dimension of the subspace of cusp forms whose $L$-functions vanish to at least a prescribed order at $s=k/2$. They determine this dimension exactly when the corresponding modular curve has genus zero (Corollary~\ref{corollary_central}).

For $N>1$, choose $q$ distinct critical indices in $\{0,\ldots,k-2\}$. By Theorem~\ref{theorem_basic} and the identification of critical values with periods in \eqref{equation_period-definition}, a sufficient condition for the corresponding period functionals to span $S_k(\Gamma_0(N))^*$ is
\[
q\ge d_k(N)+g_N-\ind_{\{k=2\}}.
\]
If $g_N=0$, the selected periods span $S_k(\Gamma_0(N))^*$ if and only if $q\ge d_k(N)$. In particular, any $d_k(N)$ distinct critical periods form a dual basis when $d_k(N)\le k-1$. Fukuhara and Yang \cite{fukuhara2009period} give explicit period bases at level two. Choie, Park, and Zagier \cite[Section~7]{choie2019periods} prove injectivity of the even and odd period maps at levels $2,3,5$ in every even weight. Theorem~\ref{theorem_basic} recovers these results and bounds the rank for arbitrary selections of critical periods, without a parity restriction.

There are only $k-1$ critical periods on the path from $0$ to $\infty$, so they cannot span the dual when $ d_k(N)>k-1$. For each fixed even $k\ge2$, this obstruction occurs for all sufficiently large $N$. We therefore allow the initial cusp to vary. For $\alpha\in\Q$, define the period integrals
\begin{equation}\label{equation_additive-period}
r_n(f,\alpha)=\int_\alpha^{i\infty}f(z)(z-\alpha)^n\,dz
             =i^{n+1}\int_0^\infty f(\alpha+iy)y^n\,dy,\qquad 0\le n\le k-2.
\end{equation}
The path is vertical, and the integral converges absolutely at both endpoints because $f$ is cuspidal. At $\alpha=0$, this agrees with \eqref{equation_period-definition}: $r_n(f,0)=r_n(f)$. We also regard periods as linear functionals of $f$, writing $r_n(\,\cdot\,,\alpha)$.

Period integrals between cusps belong to the classical theory of modular symbols and Shimura integrals \cite{manin1972parabolic,manin1973periods,shokurov1981shimura}. By Eichler--Shimura theory, the class in $H^1_{\mathrm{par}}(\Gamma_0(N),\operatorname{Sym}^{k-2}\C^2)$ determines the form \cite[Section~2]{pasol2013modular}. After splitting paths at $\infty$ and changing the polynomial basis, we conclude that the family of all $r_n(\,\cdot\,,\alpha)$, $\alpha\in\Q$, spans the dual. Finite modular symbol descriptions and explicit reconstruction from complete period data are developed
in \cite{merel1994universal} and \cite[Theorem~6.1]{pasol2013modular}. 

Related determination results use multiplicative twists. Merel \cite[Corollary~2]{merel2009symboles} determines weight-two newforms from central values of twists with conductors dividing the level, together with specified local data. Dickson and Neururer \cite[Theorems~3.2 and~3.4]{dickson2018products} give vanishing criteria using twisted critical values for elements of new spaces subject to Atkin--Lehner conditions. For level-one forms, Ni and Xue \cite[Theorems~1.1--1.4]{ni2026twisted} prove independence of twisted periods for primitive characters of odd squarefree conductor and sufficiently large weight, subject to parity and index restrictions; when the character varies, they also require distinct character values at $2$. Ni \cite[Corollary~6.1]{ni2026uniform} proves a density-one basis result for prime quadratic twists at each noncentral critical index, again for level-one forms. For primitive characters these periods are finite Gauss sum combinations of the period integrals $r_n(f,\alpha)$. 

In Section~\ref{section_periods}, we study when finite collections of the period functionals $r_n(\,\cdot\,,\alpha)$ span the dual of $S_k(\Gamma_0(N))$, allowing both the rational point $\alpha$ and the critical index $n$ to vary. Classical Hecke theory supplies a finite spanning family of the functionals $r_{k-2}(\,\cdot\,,\alpha)$ at every level in weights $k\ge4$ (Proposition~\ref{proposition_hecke-periods}). With Chebyshev systems, we can bound the rank of prescribed subfamilies on suitable real paths and select explicit bases when the relevant quotient curves have genus zero. The bounds depend on the dimensions of the relevant cusp-form spaces and the genera of these curves. For example, the periods $r_n(f,\alpha)$ at the three fixed rational points $\alpha=0,1/3,2/3$ determine every cusp form in $S_k(\Gamma_0(27))$ uniquely, in every even weight $k\ge2$. For $k\ge4$, an explicit dual basis can be selected from the corresponding period functionals (Example~\ref{remark_period-level27}). At level $35$, the periods $r_n(f,\alpha)$ at five fixed rational points likewise determine every cusp form in $S_k(\Gamma_0(35))$ uniquely in every even weight $k\ge2$ (Example~\ref{remark_period-level35}).

The paper is organized as follows. Section~\ref{section_chebyshev} develops zero bounds and interpolation rank estimates for integral transforms using Chebyshev systems. In Section~\ref{section_modular}, we apply these results to completed modular $L$-functions. We examine sharpness and derivative evaluations, and prove that the real zeros of Mellin transforms satisfy the bound of dimension minus one throughout the full nonzero cusp-form space at only finitely many levels for each fixed even weight. The section concludes with weighted moments of critical values and applications of Shimura's period relations to rational spans of eigenform periods and coefficient fields generated by period ratios. In Section~\ref{section_periods}, we study when finite families of critical periods span the dual of the cusp-form space and construct explicit dual bases. A comparison table at small levels summarizes the scope of the four geometric constructions.

\section{Chebyshev systems and integral transforms}\label{section_chebyshev}

\subsection{Basic tools}

We use the classical Chebyshev criteria from \cite{borwein1995polynomials,karlin1966tchebycheff}. 

\begin{definition}\label{definition_cheb}
A $d$-dimensional real space of continuous functions on an interval is a \emph{Chebyshev space} if every nonzero element has at most $d-1$ distinct zeros. A basis is called a \emph{Chebyshev system}, or \emph{T-system}. For functions analytic near the interval, the corresponding property with zeros counted with multiplicity is called \emph{extended Chebyshev}, or \emph{ET}. An ordered system is \emph{complete Chebyshev} (CT), respectively \emph{extended complete Chebyshev} (ECT), when every initial subsystem is T, respectively ET.
\end{definition}
A $d$-dimensional space has the ET property if and only if every Hermite interpolation problem with $d$ scalar conditions has a unique solution in that space. For linearly independent analytic functions $\phi_1,\ldots,\phi_d$ on an open interval $J$, the classical criteria take the following forms.
\begin{enumerate}
\item The system is T if and only if $\det[\phi_j(x_i)]_{i,j=1}^d\ne0$ at every choice of distinct nodes $x_1,\ldots,x_d\in J$. Applying this to each initial subsystem characterizes CT.
\item The system is ECT if and only if all its initial Wronskians satisfy
\[
\det[\phi_j^{(i-1)}(x)]_{i,j=1}^r\ne0\qquad(x\in J,\ 1\le r\le d).
\]
\item Multiplying every function by the same nowhere-zero real analytic factor, or making an analytic change of variable with nowhere-zero derivative onto another interval, preserves these properties.
\end{enumerate}

For a nonzero real analytic function $F$ with finitely many zeros on $J$, denote its number of sign changes by $S_J(F)$. Thus $S_J(F)\le Z_J(F)$. The following sign-change criterion is classical \cite{karlin1966tchebycheff}.
\begin{lemma}\label{lemma_sign-change}
Let $J$ be an open interval, and let $F$ be a nonzero real analytic function on $J$ with finitely many zeros. Suppose $\phi_1,\ldots,\phi_q$ form an ECT-system on $J$ and every $F\phi_j$ is absolutely integrable. Then
\[
\int_J F(x)\phi_j(x)\,dx=0\quad(1\le j\le q)\quad\Longrightarrow\quad S_J(F)\ge q.
\]
\end{lemma}
\begin{proof}
Suppose, to the contrary, that $F$ has $r<q$ sign changes, at $x_1<\cdots<x_r$. Choose a nonzero $P\in\Span_{\R}\{\phi_1,\ldots,\phi_{r+1}\}$ vanishing at these points. The ECT property implies that these are all the zeros of $P$ in $J$ and that they are simple. Thus $P$ changes sign exactly where $F$ does. After changing the sign of $P$ if necessary, $F P\ge0$ on $J$ and $F P>0$ on some
open subinterval. Since $F P$ is absolutely integrable, $\int_J F P>0$. This contradicts the orthogonality assumptions. 
\end{proof}

\subsection{Examples of Chebyshev systems}

For real exponential polynomials, the classical Chebyshev property \cite{karlin1966tchebycheff} holds in the following confluent form.

\begin{lemma}\label{lemma_exponential}
For real $\lambda_1<\cdots<\lambda_\ell$ and integers $m_i\ge1$,
the ordered blocks
\[
\bigl(x^a e^{\lambda_i x}\bigr)_{0\le a<m_i},\qquad 1\le i\le\ell,
\]
form an ECT-system on $\R$. Thus every nonzero member of their span has at most $\sum_i m_i-1$ real zeros, with multiplicity.
\end{lemma}

Indeed, the Wronskian of the displayed system is
\[
e^{x\sum_i m_i\lambda_i}\prod_i\prod_{a=0}^{m_i-1}a!\,\prod_{i<j}(\lambda_j-\lambda_i)^{m_i m_j},
\]
and the same formula applies to every initial subsystem. Its factors are nonzero, so the Wronskian criterion proves the assertion.

\begin{lemma}\label{lemma_kernels}
The following ordered systems are ECT-systems on $(0,\infty)$:
\begin{enumerate}
\item
$x^{a_1},\ldots,x^{a_q}$ for any real $a_1<\cdots<a_q$.

\item
$\cosh(\lambda_1x),\ldots,\cosh(\lambda_qx)$ for $0\le\lambda_1<\cdots<\lambda_q$, and
$\sinh(\lambda_1x),\ldots,\sinh(\lambda_qx)$ for $0<\lambda_1<\cdots<\lambda_q$.

\item
For $\beta\in\{0,1\}$, define
\[
K_\beta(u,x)=\sum_{j\ge0}\frac{u^j x^{\beta+2j}}{(\beta+2j)!}.
\]
For any $0\le u_1<\cdots<u_\ell$ and positive integers $m_1,\ldots,m_\ell$, consider the functions
\[
\bigl(\partial_u^aK_\beta(u_i,\cdot)\bigr)_{1\le i\le\ell,\ 0\le a<m_i},
\]
ordered first by $i$ and then by $a$.
\end{enumerate}
\end{lemma}
\begin{proof}
All zeros in this proof are counted with multiplicity. Part (1) follows from Lemma~\ref{lemma_exponential} by setting $t=\log x$. We next prove (3). Let $q=\sum_{i=1}^{\ell}m_i$. For $u_i>0$, write $z_i=\sqrt{u_i}$. By the identities
\[
K_0(u,x)=\cosh(\sqrt{u}\,x),\qquad
K_1(u,x)=\frac{\sinh(\sqrt{u}\,x)}{\sqrt{u}}\qquad(u>0)
\]
and $\partial_u=(2z)^{-1}\partial_z$, we have
\[
\partial_u^aK_\beta(u_i,x)=A_{i,a}(x)e^{z_i x}+B_{i,a}(x)e^{-z_i x},
\]
where $A_{i,a}$ and $B_{i,a}$ have degree at most $a$, and the coefficient of $x^a$ in $A_{i,a}$ is $\frac{(2z_i)^{-a}}{2z_i^\beta}\ne0.$ Differentiating the power series at $u=0$, we obtain
\[
\partial_u^aK_\beta(0,x)=\frac{a!}{(\beta+2a)!}x^{\beta+2a}\qquad(a\ge0).
\]

Suppose
\[
\sum_{i=1}^{\ell}\sum_{a=0}^{m_i-1}c_{i,a}\,\partial_u^aK_\beta(u_i,x)=0.
\]
By Lemma~\ref{lemma_exponential}, the polynomial coefficients corresponding to distinct exponents must vanish separately. For each $u_i>0$, comparison of the coefficients of $e^{z_i x}$ in decreasing order of polynomial degree shows that
$c_{i,a}=0$ for all $a$. If $u_1=0$, the remaining terms are distinct monomials, so their coefficients also vanish. Thus the functions in (3) are linearly independent.

These functions lie in the span of
\[
x^k e^{\pm z_i x},\qquad u_i>0,\quad 0\le k<m_i,
\]
together with
\[
1,x,\ldots,x^{\beta+2m_1-2}
\]
if $u_1=0$. This space has dimension at most $2q$. Hence, by Lemma~\ref{lemma_exponential}, every nonzero linear combination $G$ of the functions in (3) has at most $2q-1$ real zeros. Since $G(-x)=(-1)^\beta G(x)$, its nonzero real zeros occur in symmetric pairs with equal multiplicities, and $G(0)=0$ when $\beta=1$. Consequently,
\[
2Z_{(0,\infty)}(G)+\beta\le 2q-1,
\]
so $Z_{(0,\infty)}(G)\le q-1$. Every initial subsystem is of the same form, which proves (3).

Finally, (2) follows from (3) by taking $m_i=1$ and $u_i=\lambda_i^2$, using
\[
K_0(\lambda_i^2,x)=\cosh(\lambda_i x),\qquad
\lambda_iK_1(\lambda_i^2,x)=\sinh(\lambda_i x)\quad(\lambda_i>0),
\]
and the invariance of the ECT property under multiplication of individual functions by nonzero constants.
\end{proof}

\subsection{Integral transforms and moment maps}

In the following results, bounds of $d+\delta-1$ zeros, counted with multiplicity, for nonzero functions in a $d$-dimensional space correspond to Chebyshev properties with accuracy at most $\delta$; see \cite[Definition~1]{gavrilov2003two} and \cite{novaes2017extended}. Here we combine such bounds with the sign-change principle to obtain rank estimates for Hermite evaluation and moment maps, with the excess $\delta$ controlling the possible rank loss.

Let $J$ be an open real interval, let $I$ be a real interval with nonempty interior, and let $V$ be a $d$-dimensional real vector space of real analytic functions on $J$. Suppose every nonzero $F\in V$ has finitely many zeros on $J$.
For a real-valued kernel $K(u,x)$, define
\[
(TF)(u)=\int_J F(x)K(u,x)\,dx.
\]
Assume that, for every $F\in V$, this integral defines a function analytic on a neighborhood of $I$ and that, for every $a\ge0$,
\begin{equation}\label{equation_transfer-derivative}
(TF)^{(a)}(u)=\int_J F(x)\partial_u^aK(u,x)\,dx,\qquad u\in I,
\end{equation}
with absolute convergence, locally uniform in $u$.

\begin{proposition}\label{proposition_transfer}
Assume that, for every $u_1<\cdots<u_\ell$ in $I$ and positive integers $m_1,\ldots,m_\ell$, the functions
\[
\bigl(\partial_u^aK(u_i,\cdot)\bigr)_{1\le i\le\ell,\ 0\le a<m_i},
\]
ordered first by $i$ and then by $a$, form an ECT-system on $J$. Then $T$ is injective and
\begin{equation}\label{equation_transfer-zeros}
Z_I(TF)\le S_J(F)\qquad(0\ne F\in V).
\end{equation}
If, for some nonnegative integer $\delta$,
\[
S_J(F)\le d+\delta-1\qquad(0\ne F\in V),
\]
then every Hermite evaluation map
\[
E:V\longrightarrow\R^q,\qquad E(F)=\bigl((TF)^{(a)}(u_i)\bigr)_{1\le i\le\ell, ~0\le a<m_i},\qquad q=\sum_{i=1}^{\ell}m_i,
\]
satisfies
\begin{equation}\label{equation_transfer-rank}
\rank E\ge \min\{d,[q-\delta]_+\}.
\end{equation}
In particular, if $\delta=0$, then $T(V)$ is a $d$-dimensional ET-space on $I$.
\end{proposition}

\begin{proof}
Let $0\ne F\in V$, and suppose
\[
(TF)^{(a)}(u_i)=0\qquad(1\le i\le\ell,\ 0\le a<m_i),
\]
where $u_1<\cdots<u_\ell$ lie in $I$. The derivative identity \eqref{equation_transfer-derivative} rewrites these conditions as
\[
\int_J F(x)\partial_u^aK(u_i,x)\,dx=0\qquad(1\le i\le\ell,\ 0\le a<m_i).
\]
The kernels form an ECT-system, so Lemma~\ref{lemma_sign-change} forces $q\le S_J(F)$. If $TF\equiv0$, the choice $\ell=1$ and $m_1>S_J(F)$ would contradict this inequality. Thus $T$ is injective. Applying this inequality to any finite collection of zeros of $TF$, with their multiplicities, proves \eqref{equation_transfer-zeros}. The same argument applies to zeros at endpoints belonging to $I$, since the derivative identities and the ECT hypothesis hold there as well.

Now assume $S_J(F)\le d+\delta-1$ for every $0\ne F\in V$. The zero bound makes $E$ injective when $q\ge d+\delta$. If $q<d+\delta$, adjoin evaluations at $d+\delta-q$ additional distinct interior points. The augmented map has rank $d$, and deleting these coordinates lowers rank by at most $d+\delta-q$. Hence $\rank E\ge q-\delta$; together with nonnegativity, this proves \eqref{equation_transfer-rank}.
\end{proof}

\begin{corollary}\label{corollary_rational-sharpness}
Let $V$ be a $d$-dimensional ET-space, with $d\ge1$, of real analytic functions on a nonempty open interval $I$, and let $V_{\Q}\subset V$ be a vector space over $\Q$ whose real span is $V$. For any $d-1$ pairwise disjoint nonempty open subintervals of $I$, there is a nonzero $F\in V_{\Q}$ with exactly one simple zero in each prescribed interval and no other zeros in $I$.
\end{corollary}

\begin{proof}
For $d=1$, any nonzero element of $V_{\Q}$ suffices. Otherwise choose a point $s_j$ in each prescribed interval. The $d-1$ conditions $F_0(s_j)=0$ admit a nonzero solution in $V$. The ET property forces these to be all its zeros and to be simple. For each $j$, take $a_j<s_j<b_j$ in the same interval with $F_0(a_j)F_0(b_j)<0$. In a basis of $V$ chosen from $V_{\Q}$, approximate the coordinates of $F_0$ by rational numbers. A sufficiently close approximation $F\in V_{\Q}$ preserves these finitely many strict inequalities, and the ET bound forces precisely the required simple zeros.
\end{proof}

We next consider moments of a single vector space on several intervals. Let $V$ be a real vector space of dimension $d>0$, let $(J_j)_j$ be a finite family of open real intervals, and let $R_j:V\to C^{\omega}(J_j,\R)$ be real-linear maps. Assume that, for every $0\ne f\in V$, each $R_jf$ is nonzero and has finitely many zeros. Write $R=(R_j)_j$ and define the total number of sign changes by
\[
\mathcal S_R(f):=\sum_j S_{J_j}(R_jf)\qquad(0\ne f\in V).
\]
Choose integers $q_j\ge0$ and families $\Phi_j=(\phi_{j,a})_{a=1}^{q_j}$, required to be real ECT-systems on $J_j$ whenever $q_j>0$. With $\Phi=(\Phi_j)_j$ and $q=\sum_jq_j$, define the moment map
\[
M_\Phi:V\longrightarrow\R^q,\qquad
(M_\Phi f)_{j,a}=\int_{J_j}(R_jf)(x)\phi_{j,a}(x)\,dx,
\]
assuming absolute convergence for every $f\in V$. For an integer $D\ge q$, call this moment family \emph{extendible to $D$ moments} if $D-q$ further moments can be added on the same intervals, retaining the original coordinates, so that each nonempty kernel family forms an ECT-system after reordering and all added integrals converge absolutely for every $f\in V$.

\begin{proposition}\label{proposition_joint-rank}
For the moment map above, let $\delta\ge0$ be an integer and assume
\[
\mathcal S_R(f)\le d+\delta-1\qquad(0\ne f\in V).
\]
If $q<d+\delta$, assume that the moment family is extendible to $d+\delta$ moments. Then
\[
\min\{d,[q-\delta]_+\}\le\rank_{\R}M_\Phi\le\min(d,q).
\]
\end{proposition}
\begin{proof}
Let $0\ne f\in\ker M_\Phi$. Applying Lemma~\ref{lemma_sign-change} on each interval with $q_j>0$, we have
\[
q=\sum_jq_j\le\mathcal S_R(f)\le d+\delta-1.
\]
Thus $M_\Phi$ is injective when $q\ge d+\delta$. If $q<d+\delta$, extend the family to $d+\delta$ moments. The extended map has rank $d$, and deleting the added coordinates lowers rank by at most $d+\delta-q$. Hence $\rank_{\R}M_\Phi\ge q-\delta$. Together with $0\le\rank_{\R}M_\Phi\le\min(d,q)$, this proves the bounds.
\end{proof}
For a direct sum $V=\bigoplus_\lambda V_\lambda$, suppose that invertible row operations reduce a map $\mathcal P$ to moment maps $M_\lambda$ on the individual summands, possibly with additional zero rows. If each nonzero summand satisfies Proposition~\ref{proposition_joint-rank} with parameters $d_\lambda,q_\lambda,\delta_\lambda$, then $\rank_{\R}\mathcal P=\sum_\lambda\rank_{\R}M_\lambda$, and hence
\begin{equation}\label{equation_joint-rank}
\sum_\lambda\min\{d_\lambda,[q_\lambda-\delta_\lambda]_+\} \le \rank_{\R}\mathcal P
\le \sum_\lambda\min(d_\lambda,q_\lambda).
\end{equation}
By \eqref{equation_joint-rank}, a necessary condition for injectivity of $\mathcal P$ is $q_\lambda\ge d_\lambda$ for every $\lambda$, and a sufficient condition is $q_\lambda\ge d_\lambda+\delta_\lambda$ for every $\lambda$. If all $\delta_\lambda=0$, then $\mathcal P$ is injective if and only if $q_\lambda\ge d_\lambda$ for every $\lambda$.

The next criterion will be used for weighted moments in Section~\ref{subsection_weighted-moments}. Related criteria for preservation of Chebyshev properties under integration appear in \cite[Theorems~1.1 and~2.4]{gasull2020chebyshev} and \cite{huang2023new}.

\begin{proposition}\label{proposition_positive-measure-rank}
Let $V$ be a $d$-dimensional real space of real-valued functions analytic near $J=[a,b]$, where $a<b$ and $d\ge1$. Suppose $Z_J(F)\le d+\delta-1$ for every $0\ne F\in V$, with $\delta\in\Z_{\ge0}$, counting analytic multiplicities also at $a,b$. Let $\mu$ be a nonzero finite positive Borel measure on $J$, and write $p=\#\operatorname{supp}\mu$, allowing $p=\infty$. For a real CT-system $\phi_1,\ldots,\phi_q$ of continuous functions on the closed interval $J$, with $q\ge1$, set
\[
E_\mu(F)=\left(\int_J F(x)\phi_j(x)\,d\mu(x)\right)_{j=1}^q.
\]
Then
\[
\min\{d,[\min(q,p)-\delta]_+\}\le\rank_{\R}E_\mu\le\min(d,q,p).
\]
In particular, if $V$ is an ET-space, then
$\rank_{\R}E_\mu=\min(d,q,p)$.
\end{proposition}
\begin{proof}
The upper bound follows by factoring through evaluations at the support points when $p<\infty$.  If $p\le q$, the CT property and positivity of the atomic weights imply that the first $p$ moments are related to the $p$ support values by an invertible matrix. Thus $\ker E_\mu$ consists exactly of the functions vanishing on the support. Adjoining $[d+\delta-p]_+$ evaluations at distinct points of $J$ outside the support makes the map injective, so $\rank E_\mu\ge d-[d+\delta-p]_+$.

Suppose $p>q$. Every $F\in\ker E_\mu$ nonzero somewhere on the support has at least $q$ sign changes in $(a,b)$. Indeed, if its sign changes are $t_1<\cdots<t_s$ with $s<q$, the evaluation determinant of $\phi_1,\ldots,\phi_{s+1}$ at $t_1,\ldots,t_s,x$ defines a function $P$ in their span; for $s=0$, take $P=\phi_1$. The CT property makes $P$ vanish precisely at the $t_i$. On the connected set of increasingly ordered nodes, the evaluation determinant is continuous and nonzero, hence has constant sign. As $x$ crosses a $t_i$, the ordering permutation changes parity, so $P$ changes sign. After a sign choice, $FP\ge0$ on $J$, strictly wherever $F\ne0$. Since $FP$ is continuous and strictly positive at some point of the support, $\int_J FP\,d\mu>0$, contrary to the moment equations.

Put $r=[d+\delta-q]_+$ and adjoin the derivative evaluations
\[
F\longmapsto\left.\frac{d^jF}{dx^j}\right|_{x=a},\qquad 0\le j<r,
\]
where $j=0$ denotes the value $F(a)$. Suppose a nonzero $F$ lies in the augmented kernel. If $F$ is nonzero somewhere on the support, the preceding argument provides at least $q$ distinct zeros in $(a,b)$. Otherwise, $F$ vanishes on the entire support, which contains at least $q$ points different from $a$ because $p>q$. In either case, $F$ has at least $q$ zeros in $J\setminus\{a\}$, while the added conditions imply $\ord_a F\ge r$. Hence
\[
Z_J(F)\ge q+r\ge d+\delta,
\]
contrary to the zero bound. The augmented map is therefore injective. Removing its $r$ added coordinates lowers rank by at most $r$, so $\rank E_\mu\ge d-r=\min(d,q-\delta)$. Combining these bounds with nonnegativity proves the result.
\end{proof}

\begin{lemma}\label{lemma_complexification}
Let $V$ be a finite-dimensional real vector space. Every real-linear map $E:V\to\R^q$ has a complexification $E_{\C}$ satisfying
\[
\ker E_{\C}=(\ker E)\otimes_{\R}\C,\qquad\rank_{\C}E_{\C}=\rank_{\R}E.
\]
If $T$ embeds $V$ real-linearly into real-valued analytic functions near a real interval $I$ with nonempty interior, its complexification is also injective and, for $v=v_1+iv_2\ne0$ with $v_1,v_2\in V$,
\[
Z_I(Tv)\le Z_I(Tv_j)\qquad\text{whenever }v_j\ne0.
\]
Zeros are counted with multiplicity, including at endpoints belonging to $I$.
\end{lemma}

\begin{proof}
Separating real and imaginary parts proves the kernel identity and injectivity, and hence the rank identity. If $Tv$ and its first $r-1$ derivatives vanish at a real point, the same holds for $Tv_1$ and $Tv_2$, which proves the zero bound.
\end{proof}

\section{Applications to modular \texorpdfstring{$L$}{L}-functions}\label{section_modular}

\subsection{A geometric bound for cofinite groups}
Let $\Gamma\subset\SL_2(\R)$ be a discrete subgroup of finite covolume, and let $\chi:\Gamma\to\C^\times$ be a unitary character. Throughout this section, assume that $k\ge2$ is an integer and $S_k(\Gamma,\chi)\ne0$.  Define the weight-$k$ slash operator by
\[
(f|_k\gamma)(z)=(ad-bc)^{k/2}(cz+d)^{-k}f(\gamma z),\qquad
\gamma=\begin{pmatrix}a&b\\c&d\end{pmatrix}\in \operatorname{GL}^{+}_2(\R).
\]
Thus $f|_k\gamma=\chi(\gamma)f$ for $f\in S_k(\Gamma,\chi)$ and $\gamma\in\Gamma$. If $-I\in\Gamma$, we must have $\chi(-I)=(-1)^k$. Here $\chi=1$ denotes the trivial character of $\Gamma$.

\begin{lemma}\label{lemma_divisor}
Let $0\ne f\in S_k(\Gamma,\chi)$, and let $D$ be an effective integral divisor on $X_\Gamma$, with coefficient $D(P)$ at $P$. At each interior point $P$ with projective stabilizer of order $e$, assume
\[
eD(P)\le\ord_P f,
\]
where $\ord_P f$ denotes the vanishing order of $f$ at a lift of $P$ to $\HH$. At each cusp $c$, assume
\[
D(c)<\ord_c f,
\]
where $\ord_c f\in\R_{>0}$ is the leading exponent in the local Fourier expansion, with the parameter normalized by the cusp width. Then we have
\[
\deg D\le d_k(\Gamma,\chi)+g(\Gamma)-1-\ind_{\{k=2,\ \chi=1\}}.
\]
\end{lemma}

\begin{proof}
Let $L(D)$ be the space of meromorphic functions on $X_\Gamma$ with poles bounded by $D$. For $0\ne a\in L(D)$, its pullback to $\HH$ is $\Gamma$-invariant, so $af$ has the same transformation law as $f$.

At an interior point $P$, choose local coordinates $w$ on $\HH$ and $t=w^e$ on the quotient. A pole of order at most $D(P)$ in $t$ pulls back to a pole of order at most $eD(P)$ in $w$, which is cancelled by the zero of $f$. At a cusp $c$, the leading Fourier exponent of $af$ satisfies
\[
\ord_c(af)\ge\ord_c f-D(c)>0.
\]
Thus $af\in S_k(\Gamma,\chi)$, and multiplication by $f$ defines an injection
\[
L(D)\longrightarrow S_k(\Gamma,\chi),\qquad a\longmapsto af.
\]
By the Riemann--Roch theorem, we have
\[
\deg D+1-g(\Gamma)\le\dim_{\C}L(D)\le d_k(\Gamma,\chi).
\]
For the stronger bound when $k=2$ and $\chi=1$, let $\omega$ be the holomorphic differential on $X_\Gamma$ obtained by descending $f(z)\,dz$. At an interior point $P$ with projective stabilizer of order $e$, we have $\ord_P f=e\,\ord_P\omega+e-1$ using the local coordinate $t=w^e$. At a cusp $c$, we have $\ord_c\omega=\ord_c f-1$. Since $D$ is integral and the cusp exponents are integers in this case, the hypotheses imply
\[
D(P)\le\left\lfloor\frac{\ord_P f}{e}\right\rfloor=\ord_P\omega,
\qquad D(c)\le\ord_c f-1=\ord_c\omega.
\]
Thus $D\le\operatorname{div}(\omega)$. Since $d_2(\Gamma)=g(\Gamma)$, we conclude that
\[
\deg D\le\deg\operatorname{div}(\omega)=2g(\Gamma)-2=d_2(\Gamma)+g(\Gamma)-2.
\]
\end{proof}

At a nonelliptic interior point, the above condition is $D(P)\le\ord_P f$, and at a cusp with integral Fourier exponents, it is $D(c)\le\ord_c f-1$. These are the cases used below.

\subsection{Conditions for interpolation of Mellin transforms}\label{subsection_conditions}
To apply Lemma~\ref{lemma_divisor} along the imaginary axis, we impose the following conditions. Proposition~\ref{proposition_groups} verifies them for the groups considered below.

Retain the standing assumptions on $\Gamma$, $\chi$, and $k$. Fix $N>0$ as a normalization parameter for the Mellin transform. Assume the following:
\begin{enumerate}
\item[(a)]
Both $0$ and $\infty$ are cusps, and the Fourier expansions at these cusps have integral exponents. The integrality condition is equivalent to requiring
$\chi(\gamma_c)=\varepsilon^k$ whenever $\gamma_c$ is a primitive parabolic lift at either cusp conjugate to
\[
\varepsilon
\begin{pmatrix}
1&w\\0&1
\end{pmatrix},
\qquad \varepsilon\in\{1,-1\},\quad w>0.
\]

\item[(b)]
The space $S_k(\Gamma,\chi)$ is stable under the antilinear involution $f\mapsto f^c$, where
\[
f^c(z):=\overline{f(-\bar z)}.
\]
This involution conjugates the Fourier coefficients at infinity. A sufficient condition is $\Gamma^\#=\Gamma$ and
\[
\chi(\gamma^\#)=\overline{\chi(\gamma)}\qquad(\gamma\in\Gamma),
\]
where
\[
\begin{pmatrix}a&b\\c&d\end{pmatrix}^{\#}=
\begin{pmatrix}
a&-b\\-c&d
\end{pmatrix}.
\]
\end{enumerate}

For $z\in\HH$, write $[z]$ for its image in $\Gamma\backslash\HH\subset X_\Gamma$. For a cusp $c$, write $[c]$ for the corresponding point of $X_\Gamma$. Assume in addition that one of the following holds:
\begin{enumerate}
\item[(i)]
The map $y\mapsto[iy]$ is injective on $(0,\infty)$, and its image contains no elliptic points.  

\item[(ii)]
The Fricke matrix
\[
W_N=
\begin{pmatrix}
0&-1/\sqrt N\\
\sqrt N&0
\end{pmatrix}
\]
belongs to $\Gamma$. The map $y\mapsto[iy]$ is injective on $(N^{-1/2},\infty)$,
and its image contains no elliptic points.
\end{enumerate}

With this normalization,
\[
(f|_kW_N)(z)=N^{-k/2}z^{-k}f\bigl(-1/(Nz)\bigr).
\]
In case~(ii), $W_N$ identifies the cusps $0$ and $\infty$ and fixes the elliptic point $i/\sqrt N$, which is excluded from the divisor count.

We use the following normalization of the Mellin transform:
\begin{equation}\label{equation_mellin}
\Lambda_N(f,s)=N^{s/2}\int_0^\infty f(iy)y^{s-1}\,dy.
\end{equation}
If the cusp width at infinity is $w$ and the Fourier expansion is
\[
f(z)=\sum_{n\ge1}a_f(n)e^{2\pi inz/w},
\]
then
\begin{equation}\label{equation_width-completion}
\Lambda_N(f,s)=\left(\frac{\sqrt N\,w}{2\pi}\right)^s\Gamma(s)\sum_{n\ge1}\frac{a_f(n)}{n^s}
\end{equation}
for $\Re s$ sufficiently large.

In case~(i), changing $N$ multiplies $\Lambda_N(f,s)$ by a nonvanishing exponential in $s$, preserving zero multiplicities and the ranks of Hermite evaluation maps. In case~(ii), $N$ is determined by $W_N$.

In case~(ii), $W_N^2=-I$ and $\chi(-I)=(-1)^k$ imply $(i^k\chi(W_N))^2=1$. Define $\beta\in\{0,1\}$ by
\begin{equation}\label{equation_parity}
(-1)^\beta=i^k\chi(W_N).
\end{equation}
The corresponding functional equation is
\begin{equation}\label{equation_functional-equation}
\Lambda_N(f,s)=(-1)^\beta\Lambda_N(f,k-s).
\end{equation}
Define
\begin{equation}\label{equation_unfolded-u}
U_f(u)=\sum_{j\ge0}\frac{\Lambda_N^{(\beta+2j)}(f,k/2)}{(\beta+2j)!}\,u^j.
\end{equation}
At the center, Taylor coefficients in orders of parity opposite to $\beta$ vanish by the functional equation. With $U_f$ defined above, the Taylor series therefore factors as
\[
\Lambda_N(f,k/2+t)=t^\beta U_f(t^2).
\]

Under (a) and (b), define
\[
S_k(\Gamma,\chi)_{\R}=\{f\in S_k(\Gamma,\chi):f^c=f\}.
\]
By (b), this is a real form of $S_k(\Gamma,\chi)$, of real dimension $d_k(\Gamma,\chi)$. It consists precisely of the forms with real Fourier coefficients at infinity. These forms take real values on the imaginary axis.

We use the following integral representations to apply the Chebyshev estimates.
\begin{lemma}\label{lemma_mellin}
For every $f\in S_k(\Gamma,\chi)$, define
\[
F_f(x)=e^{kx/2}f(ie^x/\sqrt N),\qquad x\in\R.
\]
Then
\begin{equation}\label{equation_bilateral-mellin}
\Lambda_N(f,k/2+t)=\int_{\R}F_f(x)e^{tx}\,dx.
\end{equation}
The integral converges absolutely and locally uniformly for $t\in\C$, defines an entire function, and may be differentiated under the integral sign to every order. The map $f\mapsto\Lambda_N(f,\cdot)$ is injective.

In case~(ii), $F_f(-x)=(-1)^\beta F_f(x)$, and
\begin{equation}\label{equation_folded}
\Lambda_N(f,k/2+t)=
\begin{cases}
2\displaystyle\int_0^\infty F_f(x)\cosh(tx)\,dx,&\beta=0,\\[1mm]
2\displaystyle\int_0^\infty F_f(x)\sinh(tx)\,dx,&\beta=1.
\end{cases}
\end{equation}
Moreover, for $K_{\beta}(u,x)$ defined in Lemma~\ref{lemma_kernels}(3), we have
\begin{equation}\label{equation_normalized-kernel}
U_f(u)=2\int_0^\infty F_f(x)K_\beta(u,x)\,dx.
\end{equation}
This integral converges absolutely and locally uniformly for $u\in\C$ and may be differentiated under the integral sign to every order. In particular, $U_f$ is entire. The map $f\mapsto U_f$ is also injective.
\end{lemma}

\begin{proof}
From the identity
\[
F_f(-x)=i^k e^{kx/2}(f|_kW_N)(ie^x/\sqrt N)
\]
and cuspidality at $0$ and $\infty$, we deduce that for some $c>0$,
\begin{equation}\label{equation_mellin-decay}
F_f(x)=O\bigl(e^{k|x|/2}\exp(-ce^{|x|})\bigr)\qquad(|x|\longrightarrow\infty).
\end{equation}
Under the substitution $y=e^x/\sqrt N$, the Mellin integral takes the form \eqref{equation_bilateral-mellin}. The decay estimate ensures absolute and locally uniform convergence in $t\in\C$ and justifies differentiation under the integral sign to every order. Since $\chi$ is unitary and $f$ is cuspidal, $y^{k/2}|f(x+iy)|$ is $\Gamma$-invariant and bounded on $\HH$. Hence $a_f(n)=O(n^{k/2})$, and termwise integration for $\Re s>k/2+1$ proves \eqref{equation_width-completion}. The map $f\mapsto\Lambda_N(f,\cdot)$ is therefore injective by uniqueness of absolutely convergent Dirichlet series.

In case~(ii), the relation $f|_kW_N=\chi(W_N)f$ implies $F_f(-x)=(-1)^\beta F_f(x)$. Split \eqref{equation_bilateral-mellin} at zero and substitute $x\mapsto -x$ in the integral over $(-\infty,0)$. Using the parity of $F_f$, we combine the two integrals into \eqref{equation_folded}. Replacing $t$ by $-t$ then proves \eqref{equation_functional-equation}.

For $R>0$, the power series for $K_\beta$ satisfies
\[
|K_\beta(u,x)|\le(1+x)e^{\sqrt R\,x}\qquad(|u|\le R,\ x>0).
\]
Together with \eqref{equation_mellin-decay} and Cauchy's estimates on larger disks in the $u$-plane, this bound shows that the integral on the right-hand side of \eqref{equation_normalized-kernel} converges absolutely and locally uniformly in $u$ and may be differentiated under the integral sign to every order. The integral therefore defines an entire function of $u$. We expand $K_\beta$ in powers of $u$ and integrate term by term:
\[
\begin{aligned}
2\int_0^\infty F_f(x)K_\beta(u,x)\,dx
&=\sum_{j\ge0}\frac{2u^j}{(\beta+2j)!}
  \int_0^\infty F_f(x)x^{\beta+2j}\,dx\\
&=\sum_{j\ge0}
  \frac{\Lambda_N^{(\beta+2j)}(f,k/2)}{(\beta+2j)!}\,u^j
=U_f(u).
\end{aligned}
\]
For the second equality, we differentiate
\eqref{equation_bilateral-mellin} $\beta+2j$ times at $t=0$
and use the parity of $F_f$. Finally, the identity $\Lambda_N(f,k/2+t)=t^\beta U_f(t^2)$ and the injectivity of $f\mapsto\Lambda_N(f,\cdot)$ imply that $f\mapsto U_f$ is injective.
\end{proof}

For a real subspace $V\subseteq S_k(\Gamma,\chi)_{\R}$, we call
\[
\{\,s\mapsto\Lambda_N(f,s):f\in V\,\}
\]
its \emph{Mellin image}. In case~(ii), we call
\[
\{\,u\mapsto U_f(u):f\in V\,\}
\]
its \emph{normalized Mellin image}, where the normalization is given by $\Lambda_N(f,k/2+t)=t^\beta U_f(t^2)$. The functions are restricted to the real intervals under consideration. By Lemma~\ref{lemma_mellin}, both images have dimension $\dim_{\R}V$ on every real interval with nonempty interior.

\begin{theorem}\label{theorem_modular}
Let
\[
I_\beta=
\begin{cases}
(-\infty,k/2],&\beta=0,\\
(-\infty,k/2),&\beta=1.
\end{cases}
\]
Under the standing assumptions, suppose that (a), (b), and either (i) or (ii) hold. All evaluation maps below have complex rank at least
\begin{equation}\label{equation_modular-rank}
\min\bigl\{d_k(\Gamma,\chi),[q-g(\Gamma)+\ind_{\{k=2,\ \chi=1\}}]_+\bigr\},
\end{equation}
where $q$ is the total number of scalar conditions.
\begin{enumerate}
\item[(i)]
In case~(i), every $0\ne f\in S_k(\Gamma,\chi)$ satisfies
\begin{equation}\label{equation_full-zeros}
Z_{\R}(\Lambda_N(f,\cdot))+\sum_{c\in\{[0],[\infty]\}}(\ord_c f-1)\le d_k(\Gamma,\chi)+g(\Gamma)-1-\ind_{\{k=2,\ \chi=1\}}.
\end{equation}
The rank bound applies to the Hermite evaluation maps
\[
f\longmapsto\bigl(\Lambda_N^{(a)}(f,s_i)\bigr)_{1\le i\le\ell,\ 0\le a<m_i},\qquad s_1<\cdots<s_\ell,\quad q=\sum_{i=1}^{\ell}m_i,
\]
where $s_i\in\R$ and the $m_i$ are positive integers.

\item[(ii)]
In case~(ii), every $0\ne f\in S_k(\Gamma,\chi)$ satisfies
\begin{equation}\label{equation_u-zeros}
Z_{[0,\infty)}(U_f)+\ord_\infty f\le d_k(\Gamma,\chi)+g(\Gamma)-\ind_{\{k=2,\ \chi=1\}},
\end{equation}
or, equivalently,
\begin{equation}\label{equation_symmetric-total-zeros}
Z_{\R}(\Lambda_N(f,\cdot))+2\ord_\infty f\le\beta+2\bigl(d_k(\Gamma,\chi)+g(\Gamma)-\ind_{\{k=2,\ \chi=1\}}
 \bigr).
\end{equation}
The rank bound applies to the Hermite evaluation maps
\begin{equation}\label{equation_u-hermite}
f\longmapsto\bigl(U_f^{(a)}(u_i)\bigr)_{1\le i\le\ell,\ 0\le a<m_i},\qquad 0\le u_1<\cdots<u_\ell,\quad q=\sum_{i=1}^{\ell}m_i,
\end{equation}
where the $m_i$ are positive integers.
Equivalently, the same bound holds for
\[
f\longmapsto\left(
\bigl(\Lambda_N^{(\beta+2j)}(f,k/2)\bigr)_{0\le j<m_0},
\bigl(\Lambda_N^{(a)}(f,s_i)\bigr)_{1\le i\le\ell,\ 0\le a<m_i}
\right),
\]
where $m_0\ge0$, $m_i\ge1$, the noncentral nodes $s_i$ have pairwise distinct distances from $k/2$, and $q=m_0+\sum_i m_i$.

It also applies to
\begin{equation}\label{equation_symmetry-rank}
f\longmapsto\bigl(\Lambda_N(f,s_1),\ldots,\Lambda_N(f,s_q)\bigr),\qquad s_j\in\R,
\end{equation}
provided that the numbers $|s_j-k/2|$ are distinct and $s_j\ne k/2$ when $\beta=1$. Any distinct nodes in $I_\beta$ satisfy these conditions.
\end{enumerate}

In either case, a sufficient condition for injectivity of each of these maps is
\[
q\ge d_k(\Gamma,\chi)+g(\Gamma)-\ind_{\{k=2,\ \chi=1\}}.
\]
If $g(\Gamma)=0$, the Mellin image of $S_k(\Gamma,\chi)_{\R}$ is ET on $\R$ in case~(i). In case~(ii), the normalized Mellin image is ET on $[0,\infty)$, and the Mellin image is ET on $(-\infty,k/2)$.
\end{theorem}

\begin{proof}
For simplicity, write
\[
d=d_k(\Gamma,\chi),\qquad\delta=g(\Gamma)-\ind_{\{k=2,\ \chi=1\}}.
\]
By our assumptions, $\delta\ge0$. For $J=\R$ or $(0,\infty)$, the map $f\mapsto F_f|_J$ is injective. It therefore identifies $S_k(\Gamma,\chi)_{\R}$ with a $d$-dimensional real vector space of real analytic functions on $J$. The change of variable $y=e^x/\sqrt N$ and the nonvanishing factor $e^{kx/2}$ preserve zero multiplicities.

First let $0\ne f\in S_k(\Gamma,\chi)_{\R}$. In case~(i), choose any finite collection of zeros $iy_j$ of $f$ on the positive imaginary axis and apply Lemma~\ref{lemma_divisor} to
\[
D=\sum_j\ord_{iy_j}(f)[iy_j]+\sum_{c\in\{[0],[\infty]\}}(\ord_c f-1)[c].
\]
The interior points are distinct on $X_\Gamma$ and nonelliptic by (i). Condition~(a) and the cuspidality of $f$ imply that $\ord_c f-1$ is a nonnegative integer for each $c\in\{[0],[\infty]\}$, and each distinct cusp is counted only once. Since the collection of zeros was arbitrary, we obtain
\begin{equation}\label{equation_geometric-full-budget}
S_{\R}(F_f)\le Z_{\R}(F_f)\le d+\delta-1-\sum_{c\in\{[0],[\infty]\}}(\ord_c f-1).
\end{equation}
Thus $F_f$ has finitely many zeros and $S_{\R}(F_f)\le d+\delta-1$. In \eqref{equation_bilateral-mellin}, the functions $x^a e^{t_i x}$ form an ECT-system on $\R$ by Lemma~\ref{lemma_exponential}, while absolute convergence and differentiation under the integral sign were established in Lemma~\ref{lemma_mellin}. Proposition~\ref{proposition_transfer} and \eqref{equation_geometric-full-budget} therefore imply
\[
Z_{\R}(\Lambda_N(f,\cdot))\le S_{\R}(F_f)\le d+\delta-1-\sum_{c\in\{[0],[\infty]\}}(\ord_c f-1),
\]
which is \eqref{equation_full-zeros}. For the rank estimate, we use the uniform bound $S_{\R}(F_f)\le d+\delta-1$ in Proposition~\ref{proposition_transfer}. Every Hermite evaluation map thus has real rank at least $\min\{d,[q-\delta]_+\}$, as asserted in \eqref{equation_modular-rank} on $S_k(\Gamma,\chi)_{\R}$.

In case~(ii), choose a finite collection of zeros $iy_j$ with $y_j>N^{-1/2}$ and use
\[
D=\sum_j\ord_{iy_j}(f)[iy_j]+(\ord_\infty f-1)[\infty].
\]
By Lemma~\ref{lemma_divisor}, we have
\begin{equation}\label{equation_geometric-half-budget}
S_{(0,\infty)}(F_f)\le Z_{(0,\infty)}(F_f)\le d+\delta-\ord_\infty f.
\end{equation}
Hence $F_f$ has finitely many zeros on $(0,\infty)$ and, since $\ord_\infty f\ge1$,
\[
S_{(0,\infty)}(F_f)\le d+\delta-1.
\]
Similarly, we apply Proposition~\ref{proposition_transfer} to \eqref{equation_normalized-kernel}, using Lemma~\ref{lemma_kernels}(3) and Lemma~\ref{lemma_mellin}, to obtain
\[
Z_{[0,\infty)}(U_f)\le S_{(0,\infty)}(F_f)\le d+\delta-\ord_\infty f.
\]
The rank bound follows as above from $S_{(0,\infty)}(F_f)\le d+\delta-1$, with evaluations at $u=0$ allowed by Lemma~\ref{lemma_kernels}(3).

For the equivalent mixed map, put $t_i=s_i-k/2$ and $u_i=t_i^2$. Comparing Taylor coefficients in $\Lambda_N(f,k/2+t)=t^\beta U_f(t^2)$, we obtain
\[
\Lambda_N^{(\beta+2j)}(f,k/2)=\frac{(\beta+2j)!}{j!}U_f^{(j)}(0).
\]
At a noncentral node, $\Lambda_N^{(a)}(f,s_i)$ equals $2^a t_i^{\beta+a}U_f^{(a)}(u_i)$ plus a linear combination of lower derivatives of $U_f$ at $u_i$. Since $t_i\ne0$, these are invertible triangular changes of coordinates from the Hermite map for $U_f$. The identity
\[
\Lambda_N(f,k/2+t)=t^\beta U_f(t^2)
\]
shows that a zero of $U_f$ at $u_0>0$ corresponds to zeros of the same multiplicity at $t=\pm\sqrt{u_0}$. At the center, we have
\begin{equation}\label{equation_central-order}
\ord_{s=k/2}\Lambda_N(f,s)=\beta+2\ord_{u=0}U_f.
\end{equation}
Counting these zeros with multiplicity, we obtain
\begin{equation}\label{equation_zero-identity}
Z_{\R}(\Lambda_N(f,\cdot))=\beta+2Z_{[0,\infty)}(U_f),
\end{equation}
and this proves \eqref{equation_symmetric-total-zeros}. The map $s\mapsto(s-k/2)^2$ is injective on $I_\beta$, so the corresponding nodes $u_j=(s_j-k/2)^2$ are distinct. When $\beta=0$, we have $\Lambda_N(f,s_j)=U_f(u_j)$; when $\beta=1$, we have $\Lambda_N(f,s_j)=(s_j-k/2)U_f(u_j)$ with $s_j-k/2\ne0$. Thus the two evaluation maps have the same rank. All evaluation maps considered above take real values on $S_k(\Gamma,\chi)_{\R}$. Lemma~\ref{lemma_complexification} therefore extends their rank bounds to $S_k(\Gamma,\chi)$.

For the zero bounds, choose a nonzero real component $f_j$ of $f$. Conjugation preserves the order at infinity. At the cusp $0$, the matrix $S=\left(\begin{smallmatrix}0&-1\\1&0\end{smallmatrix}\right)$ satisfies $S^\#=-S$, and hence
\[
f^c|_kS=(-1)^k(f|_kS)^c.
\]
It therefore preserves the leading Fourier exponent at $0$ as well. Consequently,
\[
\ord_c f_j\ge\ord_c f\qquad(c=0,\infty).
\]
Lemma~\ref{lemma_complexification} bounds the zeros of $\Lambda_N(f,\cdot)$, and of $U_f$ in case~(ii), by those of the corresponding transform of $f_j$. Applying the real bounds to $f_j$ and using these cusp-order inequalities proves \eqref{equation_full-zeros} and \eqref{equation_u-zeros} for complex $f$. Finally, \eqref{equation_symmetric-total-zeros} follows from \eqref{equation_u-zeros} through the identity \eqref{equation_zero-identity}.

When $g(\Gamma)=0$, it follows that $\delta=0$. The cusp orders are at least one, so \eqref{equation_full-zeros} and \eqref{equation_u-zeros} bound the relevant zero counts by $d-1$. The injectivity in Lemma~\ref{lemma_mellin} then proves the ET assertions for the Mellin image in case~(i) and the normalized Mellin image in case~(ii). In case~(ii), the identity $\Lambda_N(f,k/2+t)=t^\beta U_f(t^2)$ transfers the same bound to $t<0$, where the change of variable and the nonvanishing factor preserve zero multiplicities.
\end{proof}

\begin{remark}\label{remark_mixed-rank}
Under the hypotheses of Theorem~\ref{theorem_modular}, write
$d=d_k(\Gamma,\chi)$ and $\delta=g(\Gamma)-\ind_{\{k=2,\ \chi=1\}}$.
For an integer $p\ge0$, define $A_p(f)=(a_f(1),\ldots,a_f(p))$, using the Fourier expansion at infinity in \eqref{equation_width-completion}.
If $E$ is any evaluation map with $q$ coordinates in Theorem~\ref{theorem_modular} or, in case~(ii), Theorem~\ref{theorem_derivatives}, then
\[
\rank\bigl(f\longmapsto(A_p(f),E(f))\bigr)
\ge\min\{d,[p+q-\delta]_+\}.
\]
Indeed, $A_p(f)=0$ implies $\ord_\infty f\ge p+1$. For a nonzero real form in the joint kernel, \eqref{equation_geometric-full-budget} or \eqref{equation_geometric-half-budget}, together with the Chebyshev sign-change argument, gives $p+q\le d+\delta-1$. Adding admissible evaluation conditions and applying Lemma~\ref{lemma_complexification} yields the stated bound. When $\delta=0$, the rank is exactly $\min(d,p+q)$.
\end{remark}

\subsection{Congruence groups and Atkin--Lehner quotients}

In this section, we study the standard congruence subgroups $\Gamma_0(N)$, $\Gamma_1(N)$, and $\Gamma(N)$. Let $B_N$ denote the Atkin--Lehner group of $X_0(N)$. For an exact divisor $e\parallel N$, the corresponding coset consists of matrices
\begin{equation}\label{equation_atkin-matrix}
\frac1{\sqrt e}
\begin{pmatrix}ae&b\\cN&de\end{pmatrix},
\qquad a,b,c,d\in\Z, \qquad ade-\frac{bcN}{e}=1.
\end{equation}
On forms of even weight and trivial nebentypus, these cosets act as commuting involutions; see \cite{atkin1970hecke}. For any subgroup $A\le B_N$, let $\Gamma_A$ be the inverse image of $A$ in the Atkin--Lehner extension of
$\Gamma_0(N)$. We regard a character $\chi:A\to\{\pm1\}$ as a character of $\Gamma_A$ via the quotient map $\Gamma_A\to A$. Then
\begin{equation}\label{equation_atkin-setting}
S_k(\Gamma_A,\chi)=S_k(\Gamma_0(N))^{A,\chi},\qquad X_{\Gamma_A}=X_0(N)/A,
\end{equation}
where $S_k(\Gamma_0(N))^{A,\chi}$ denotes the full simultaneous $\chi$-eigenspace.
Dimension formulas for $A=\langle W_N\rangle$, $N>1$, are given in \cite{zhang2022dimension}. For $N=1$, this construction gives $\Gamma_A=\SL_2(\Z)$ and the trivial character.

\begin{proposition}\label{proposition_groups}
For each pair $(\Gamma,\chi)$ in the table with $d_k(\Gamma,\chi)>0$, conditions~(a) and~(b) hold. The first three rows satisfy the hypotheses of case~(i) when $N>1$ and case~(ii) when $N=1$. In the last row, case~(i) holds when $W_N\notin A$, and case~(ii) holds when $W_N\in A$. Thus Theorem~\ref{theorem_modular} applies to all these spaces.
\begin{center}
\begin{tabular}{llll}
$\Gamma$ & $\chi$ & $X_\Gamma$ & Width at $\infty$\\ \hline
$\Gamma_0(N)$ & $1$ & $X_0(N)$ & $1$\\
$\Gamma_1(N)$ & $1$ & $X_1(N)$ & $1$\\
$\Gamma(N)$ & $1$ & $X(N)$ & $N$\\
$\Gamma_A$ & $\chi:A\to\{\pm1\}$ & $X_0(N)/A$ & $1$
\end{tabular}
\end{center}
\end{proposition}

\begin{proof}
The first three groups have finite index in $\SL_2(\Z)$. The cusps $0$ and $\infty$ are regular, with width $N$ at $0$ and widths $1,1,N$, respectively, at $\infty$; see \cite[Sections~6.3.3--6.3.4]{cohen2017modular}. Since the characters are trivial, the local Fourier multipliers are trivial, which verifies~(a). Reflection preserves the defining congruences, and the characters are trivial, so~(b) also holds. Let $S=\begin{pmatrix}0&-1\\1&0\end{pmatrix}$.
The standard fundamental domain \cite[Theorem~4.3.2]{cohen2017modular} shows that an element $\gamma\in\SL_2(\Z)$ satisfying $\gamma(iy)=iy'$ for some $y,y'>0$ must belong to $\{\pm I,\pm S\}$. For $N>1$, neither $S$ nor $-S$ belongs to $\Gamma_0(N)$. Thus, for each of the first three groups, the positive imaginary axis projects injectively and its image contains no elliptic
points, proving~(i). For $N=1$, all three groups equal $\SL_2(\Z)$ and $W_1=S$, and the ray $\{iy:y>1\}$ lies in the interior of the standard fundamental domain, proving~(ii).

The group $\Gamma_A$ is a finite union of cosets of $\Gamma_0(N)$ of the form \eqref{equation_atkin-matrix}; hence it is discrete and cofinite. Its compactified quotient is $X_0(N)/A$, and the transformation laws give \eqref{equation_atkin-setting}. Cuspidality is unchanged under this finite-index extension, since the corresponding local parameters differ by integral powers. Reflection preserves each Atkin--Lehner coset by replacing $b,c$ with $-b,-c$, so it acts trivially on $A$. Since $\chi$ is real, condition~(b) follows.

If an element of $\Gamma_A$ fixes $\infty$ or $0$, then $c=0$ or $b=0$, respectively, in \eqref{equation_atkin-matrix}. In either case, the determinant condition gives $ade=1$, so $e=1$ and the element belongs to $\Gamma_0(N)$. Thus the stabilizers of both cusps are unchanged, the width at $\infty$ is one, and $\chi$ is trivial on both stabilizers. This proves~(a).

We next check the arc conditions for $\Gamma_A$. Let $\gamma$ be a matrix of the form \eqref{equation_atkin-matrix}, and suppose that $\gamma(iy)$ lies on the positive imaginary axis for some $y>0$. Its real part vanishes, so
\begin{equation}\label{equation_axis-equation}
acNy^2+bd=0.
\end{equation}
We claim that $\gamma=\pm I$ or $\gamma=\pm W_N$. Suppose first that $a,b,c,d$ are all nonzero. By \eqref{equation_axis-equation}, we have
\[
(ade)\left(\frac{bcN}{e}\right) =(acN)(bd) =-(acN)^2y^2<0.
\]
Thus $ade$ and $bcN/e$ are nonzero integers of opposite signs, which is incompatible with the determinant condition $ade-bcN/e=1$. Hence at least one of $a,b,c,d$ vanishes. Equation~\eqref{equation_axis-equation} then implies $ac=bd=0$. Since $\gamma$ is invertible, this leaves only $b=c=0$ or $a=d=0$. In the first case, $ade=1$, so $e=1$ and $a=d=\pm1$; hence $\gamma=\pm I$. In the second case, $-bc(N/e)=1$. Since $N/e$ is a positive integer, we obtain $e=N$ and $bc=-1$, so $\gamma=\pm W_N$.

If $W_N\notin A$, only $\pm I$ can identify points on the positive imaginary axis. The axis therefore projects injectively to $X_{\Gamma_A}$, and its image contains no elliptic points, proving~(i). If $W_N\in A$, then $W_N(iy)=i/(Ny)$ exchanges the arcs $0<y<N^{-1/2}$ and $y>N^{-1/2}$, and fixes only $y=N^{-1/2}$. Hence the ray $y>N^{-1/2}$ projects injectively to $X_{\Gamma_A}$, and its image contains no elliptic points, proving~(ii).
\end{proof}

Theorem~\ref{theorem_basic} follows from the first row of Proposition~\ref{proposition_groups} and Theorem~\ref{theorem_modular}. For $N>1$, use case~(i) and the nonnegativity of the cusp contributions. For $N=1$, use case~(ii) and $\Lambda_N(f,k/2+t)=t^\beta U_f(t^2)$: on $t<0$, the change of variable $t\mapsto t^2$ and the nonzero factor $t^\beta$ preserve zero multiplicities and evaluation ranks.

Since $S_k(\Gamma_0(N))$ has a basis with rational Fourier coefficients, its Mellin image satisfies the hypotheses of Corollary~\ref{corollary_rational-sharpness} whenever it is ET.

\begin{example}\label{proposition_hyperelliptic-mellin}
Let
\[
(N,e)\in\{(22,11),(28,7),(30,15),(33,11),(46,23)\}.
\]
By \cite[Theorem~2, p.~451]{ogg1974hyperelliptic}, $W_e$ is the hyperelliptic involution of $X_0(N)$, so
\[
g\bigl(X_0(N)/\langle W_e\rangle\bigr)=0.
\]
Set $A=\langle W_e\rangle$. Since $1<e<N$, we have $W_N\notin A$. Proposition~\ref{proposition_groups} verifies the hypotheses of case~(i) for both characters of $A$. Thus, for every even weight $k\ge2$ and $\varepsilon\in\{\pm1\}$, Theorem~\ref{theorem_modular}(i) shows that the real eigenspace
\[
\{f\in S_k(\Gamma_0(N))_{\R}:f|_kW_e=\varepsilon f\},
\]
when nonzero, has an ET Mellin image on $\R$.

In weight two, the $W_e$-invariant holomorphic differentials descend to the genus-zero quotient, so the $+1$ eigenspace is zero and the $-1$ eigenspace is all of $S_2(\Gamma_0(N))$. Its Mellin image is therefore a $g_N$-dimensional ET-space on $\R$. 
\end{example}

\subsection{Fricke quotient}\label{subsection_sharpness}
We now examine the sharpness of the general zero bound in positive genus. Write
\[
h_N=g\bigl(X_0(N)/\langle W_N\rangle\bigr)
\]
for the genus of the Fricke quotient.

\begin{theorem}\label{theorem_fricke-sharpness}
Let $N>4$ and let $k\ge2$ be even, with $d_k(N)>0$. Then
\begin{equation}\label{equation_fricke-sharpness}
\begin{aligned}
\max\{d_k(N)-1,\,d_k(N)+g_N-1-\ind_{\{k=2\}}-2h_N\}
&\le \max_{0\ne f\in S_k(\Gamma_0(N))_{\R}}
Z_{\R}(\Lambda_N(f,\cdot))\\
&\le d_k(N)+g_N-1-\ind_{\{k=2\}}.
\end{aligned}
\end{equation}
If $h_N=0$, the maximum equals the upper bound and is attained by a form with rational Fourier coefficients and only simple real zeros.
\end{theorem}

\begin{proof}
The upper bound follows from Theorem~\ref{theorem_modular}(i). Vanishing at any $d_k(N)-1$ distinct real nodes imposes $d_k(N)-1$ homogeneous conditions on a $d_k(N)$-dimensional real space. These conditions have a nonzero solution, proving the first lower bound $d_k(N)-1$.

By \cite[Proposition~3.1]{martin2023root}, the expression
\[
(-1)^{k/2}\operatorname{tr}\bigl(W_N\mid S_k(\Gamma_0(N))\bigr)
+\ind_{\{k=2\}}
\]
is independent of the even weight $k\ge2$. In weight two, the invariant holomorphic differentials descend to the Fricke quotient, so the two eigenspaces have dimensions $h_N$ and $g_N-h_N$. The expression therefore equals $g_N+1-2h_N$ at $k=2$. Since it is independent of $k$, the trace satisfies
\begin{equation}\label{equation_fricke-trace}
\operatorname{tr}\bigl(W_N\mid S_k(\Gamma_0(N))\bigr)
=(-1)^{k/2}(g_N+1-2h_N)+\ind_{\{k=2\}}.
\end{equation}
Let $d_+$ be the dimension of the eigenspace with eigenvalue $(-1)^{k/2}$, whose Mellin transforms have functional-equation sign $+1$. The trace identity \eqref{equation_fricke-trace} then implies
\[
2d_+-2=d_k(N)+g_N-1-\ind_{\{k=2\}}-2h_N.
\]
The quotient map $X_0(N)\to X_0(N)/\langle W_N\rangle$ is ramified at the point represented by $i/\sqrt N$. By the Riemann--Hurwitz formula, its ramification divisor has degree $2(g_N+1-2h_N)>0$, so $g_N+1-2h_N\ge1$. The preceding identity then gives $2d_+\ge d_k(N)>0$, hence $d_+\ge1$. This eigenspace is preserved by $f\mapsto f^c$. For any distinct positive numbers $t_1,\ldots,t_{d_+-1}$, there is a nonzero real form satisfying
\[
\Lambda_N(f,k/2+t_j)=0\qquad(1\le j\le d_+-1).
\]
Its functional equation supplies the zeros at $k/2-t_j$, proving the other lower bound.

If $h_N=0$, these distinct symmetric zeros attain the upper bound, so they are simple and there are no others. Choose disjoint intervals around these zeros so that $\Lambda_N(f,\cdot)$ has opposite signs at the two endpoints of each interval. Approximating in a rational Fourier basis preserves these finitely many strict sign inequalities. The same upper bound forces exactly one simple zero in each interval and no others. If the bound is zero, any nonzero rational form suffices.
\end{proof}

When $h_N=0$, it follows that $d_k(N)+g_N-\ind_{\{k=2\}}$ is the least number of distinct real evaluations that always determine a form. If $X_0(N)$ has genus one, its nontrivial Fricke involution has a fixed point, so the quotient $X_0(N)/\langle W_N\rangle$ has genus zero by the Riemann--Hurwitz formula. At these levels, the maximal number of real zeros is therefore zero in weight two and $d_k(N)$ in every even weight $k\ge4$. 

\begin{corollary}\label{corollary_et-finiteness}
For each fixed even weight $k\geq 2$, there are only finitely many levels $N>1$ for which the full space $S_k(\Gamma_0(N))_{\R}$ is nonzero and its Mellin image is an ET-space on $\R$.
\end{corollary}

\begin{proof}
It suffices to consider $N>4$. Let $R_N$ be the number of fixed points of $W_N$ on $X_0(N)$. Each fixed point contributes one to the ramification divisor of the degree-two quotient map $X_0(N)\mapsto X_0(N)/\langle W_N\rangle$. Thus $R_N=2g_N+2-4h_N$ by the Riemann--Hurwitz formula. If the Mellin image is ET, its zero bound is $d_k(N)-1$, so Theorem~\ref{theorem_fricke-sharpness} implies
\[
R_N\le2+2\ind_{\{k=2\}}\le4.
\]
For a negative quadratic discriminant $D$, write $h(D)$ for the class number of the imaginary quadratic order of discriminant $D$. The Fricke fixed-point formula \cite[p.~454]{ogg1974hyperelliptic} expresses $R_N$ in terms of these class numbers:
\[
R_N=h(-4N)+\ind_{\{N\equiv3\pmod4\}}h(-N).
\]
In particular, $h(-4N)\le4$. Since the class number of an order is at least that of its field, only finitely many underlying imaginary quadratic fields can occur (c.f. \cite{watkins2004class}). For each such field, the class-number formula for orders bounds the conductor when the class number is at most four. Hence only finitely many levels $N$ can occur.
\end{proof}

\subsection{Rank bounds for derivative evaluations}\label{section_derivatives}

Throughout this subsection, assume the hypotheses of Theorem~\ref{theorem_modular}(ii). Let $\beta$ be defined by \eqref{equation_parity}.
In particular, the results apply to the simultaneous Atkin--Lehner eigenspaces in Proposition~\ref{proposition_groups} when $W_N\in A$.

\begin{theorem}\label{theorem_derivatives}
Each of the following maps $S_k(\Gamma,\chi)\to\C^q$ has rank at least
\[
\min\{d_k(\Gamma,\chi),[q-g(\Gamma)+\ind_{\{k=2,\ \chi=1\}}]_+\}.
\]
\begin{enumerate}
\item
For $s_0\in\R$ and integers $0\le j_1<\cdots<j_q$,
\[
f\longmapsto\bigl(\Lambda_N^{(\beta+2j_\nu)}(f,s_0)\bigr)_{\nu=1}^q.
\]
If $s_0\ne k/2$, the same bound holds with the orders
$\beta+2j_\nu$ replaced by $1-\beta+2j_\nu$.

\item
For an integer $r\ge0$ and real points $s_1,\ldots,s_q$ such that $|s_1-k/2|,\ldots,|s_q-k/2|$ are pairwise distinct, excluding $s_\nu=k/2$ when $\beta+r$ is odd,
\[
f\longmapsto\bigl(\Lambda_N^{(r)}(f,s_\nu)\bigr)_{\nu=1}^q.
\]

\end{enumerate}
If $g(\Gamma)=\ind_{\{k=2,\ \chi=1\}}$ and $q=d_k(\Gamma,\chi)$, each of these maps is an isomorphism.
\end{theorem}

\begin{proof}
For simplicity, write
\[
d=d_k(\Gamma,\chi),\qquad\delta=g(\Gamma)-\ind_{\{k=2,\ \chi=1\}}\ge0.
\]
For parts~(1) and~(2), we first work on $S_k(\Gamma,\chi)_{\R}$. By the zero bound \eqref{equation_geometric-half-budget}, the real-valued
function $F_f\not\equiv0$ associated with each $f\ne0$ has finitely many zeros and satisfies
\[
S_{(0,\infty)}(F_f)\le d+\delta-1.
\]
We apply Proposition~\ref{proposition_joint-rank} over $\R$ on the single interval $(0,\infty)$. All the required integrals converge absolutely, and differentiation under the integral is justified by Lemma~\ref{lemma_mellin}.

For part~(1), we differentiate \eqref{equation_folded} and set $t=s_0-k/2$:
\[
\Lambda_N^{(\beta+2j)}(f,s_0)=2\int_0^\infty F_f(x)x^{\beta+2j}\cosh((s_0-k/2)x)\,dx.
\]
The kernels form an ECT-system by Lemma~\ref{lemma_kernels}(1), since the common factor $\cosh((s_0-k/2)x)$ is nonvanishing on $(0,\infty)$. For the other parity, the kernels are
\[
x^{1-\beta+2j}\sinh((s_0-k/2)x).
\]
When $s_0\ne k/2$, the common hyperbolic factor is again nonvanishing on $(0,\infty)$, so the same argument applies.

For part~(2), let $\lambda=|s-k/2|$. Up to nonzero constant factors, the kernels are
\[
\begin{cases}
x^r\cosh(\lambda x),&\beta+r\text{ even},\\
x^r\sinh(\lambda x),&\beta+r\text{ odd}.
\end{cases}
\]
In the odd case, $\lambda=0$ is excluded because the kernel vanishes identically.
After ordering the distinct values of $\lambda$, Lemma~\ref{lemma_kernels}(2) applies. Multiplication by $x^r$ preserves the ECT property on $(0,\infty)$, and the nonzero constant factors do not affect rank. In both parts, adding unused orders or distinct positive distances extends the ECT-system to any larger size. We may therefore apply Proposition~\ref{proposition_joint-rank} over $\R$ and then pass to $\C$ by Lemma~\ref{lemma_complexification}.

\end{proof}

We now filter the cusp-form space by the vanishing order of its Mellin transforms at the center.

\begin{corollary}\label{corollary_central}
Write $d=d_k(\Gamma,\chi)$ and $\delta=g(\Gamma)-\ind_{\{k=2,\ \chi=1\}}$. For an integer $r\ge0$, put
\[
\mathcal F_r=\left\{f\in S_k(\Gamma,\chi):
\Lambda_N^{(\beta+2j)}(f,k/2)=0\text{ for }0\le j<r\right\}.
\]
Then
\begin{equation}\label{equation_central-filtration}
[d-r]_+\le\dim\mathcal F_r\le\min\{d,[d+\delta-r]_+\}.
\end{equation}
In particular, $\dim\mathcal F_r=[d-r]_+$ when $\delta=0$.
For $f\ne0$, membership in $\mathcal F_r$ is equivalent to
$\ord_{s=k/2}L(f,s)\ge\beta+2r$.
\end{corollary}
\begin{proof}
The dimension bounds follow from the rank--nullity theorem and Theorem~\ref{theorem_derivatives}(1). The functional equation forces the central derivatives of the opposite parity to vanish, and the completion factor is holomorphic and nonzero at $k/2$, proving the equivalence.
\end{proof}

\subsection{Weighted moments of critical values}\label{subsection_weighted-moments}
We now illustrate how the integral-transform principle of Proposition~\ref{proposition_positive-measure-rank} can be used to study the linear independence of other families of linear functionals on spaces of cusp forms. We apply it to the weighted moments of critical $L$-values introduced by Ono, Rolen, and Sprung \cite{ono2017zeta}. In this case, the functional equation expresses suitable centered moments in terms of $U_f$ and a positive discrete measure, so the zero bounds for $U_f$ yield the required rank estimates.

Let $k\ge2$ be even, and let $S_k(\Gamma,\chi)\ne0$ satisfy the hypotheses of Theorem~\ref{theorem_modular}(ii). Retain the sign $(-1)^\beta$ of the functional equation from \eqref{equation_parity}. Following \cite[(1.1)]{ono2017zeta}, for integers $m\ge0$ we define
\[
M_f(m)=\frac1{(k-2)!}\sum_{j=0}^{k-2}\binom{k-2}{j}\Lambda_N(f,j+1)j^m.
\]
Here and below, we use the convention $0^0=1$.

\begin{proposition}\label{proposition_weighted-moment-rank}
For every integer $q\ge1$, the map on $S_k(\Gamma,\chi)$ satisfies
\begin{equation}\label{equation_weighted-moment-rank}
\begin{aligned}
\min\{d_k(\Gamma,\chi),[\min(q,k/2-\beta)-g(\Gamma)]_+\}
&\le\rank\bigl(f\longmapsto(M_f(\beta+2r))_{r=0}^{q-1}\bigr)\\
&\le\min\{d_k(\Gamma,\chi),q,k/2-\beta\}.
\end{aligned}
\end{equation}
A necessary condition for injectivity is $\min(q,k/2-\beta)\ge d_k(\Gamma,\chi)$, and a sufficient condition is $\min(q,k/2-\beta)\ge d_k(\Gamma,\chi)+g(\Gamma)$. In genus zero the rank equals $\min\{d_k(\Gamma,\chi),q,k/2-\beta\}$; the functionals
$f\mapsto M_f(\beta+2r)$, $0\le r<d_k(\Gamma,\chi)$, form a dual basis if and only if $d_k(\Gamma,\chi)\le k/2-\beta$.
\end{proposition}
\begin{proof}
For $k=2$, we have $M_f(0)=\Lambda_N(f,1)$ and $M_f(m)=0$ for $m\ge1$. If $\beta=1$, the selected moment map is zero, as are both rank bounds. If $\beta=0$, it has the same rank as the central evaluation $f\mapsto\Lambda_N(f,1)$;
the lower bound follows from Theorem~\ref{theorem_modular}(ii), and the upper bound is immediate.

Henceforth we assume $k\ge4$. Set
\[
C_f(m)=\frac1{(k-2)!}\sum_{j=0}^{k-2}\binom{k-2}{j}\Lambda_N(f,j+1)(j+1-k/2)^m.
\]
By the functional equation, $C_f(m)=0$ for
$m\not\equiv\beta\pmod2$. After pairing the terms with
$j+1-k/2=\pm\ell$, we can write the centered moments as
\[
C_f(\beta+2r)=\int_{[0,(k/2-1)^2]}x^rU_f(x)\,d\nu_\beta(x),\qquad r\ge0,
\]
where $\delta_x$ denotes unit mass at $x$ and
\[
\nu_\beta=\frac1{(k-2)!}\left((1-\beta)\binom{k-2}{k/2-1}\delta_0+2\sum_{\ell=1}^{k/2-1}\binom{k-2}{k/2-1+\ell}\ell^{2\beta}\delta_{\ell^2}\right).
\]
Since $k\ge4$, the interval is nondegenerate and the support of $\nu_\beta$ has $k/2-\beta\ge1$ points, all with positive mass. By Lemma~\ref{lemma_mellin}, the normalized Mellin image of $S_k(\Gamma,\chi)_{\R}$ has dimension $d_k(\Gamma,\chi)$. Every nonzero function $U_f$ in this image satisfies $Z_{[0,(k/2-1)^2]}(U_f)\le d_k(\Gamma,\chi)+g(\Gamma)-1$ by Theorem~\ref{theorem_modular}(ii), since $\ord_\infty f\ge1$. Applying Proposition~\ref{proposition_positive-measure-rank} to this image with kernels $1,x,\ldots,x^{q-1}$ and then complexifying by Lemma~\ref{lemma_complexification} proves the asserted bounds for the centered moments. Expanding powers of $j=(j+1-k/2)+(k/2-1)$ by the binomial theorem, we can express the original moments in terms of the centered moments:
\begin{equation}\label{equation_moment-binomial-change}
M_f(m)=\sum_{a=0}^m\binom{m}{a}(k/2-1)^{m-a}C_f(a).
\end{equation}
Since $C_f(a)=0$ unless $a\equiv\beta\pmod2$, we obtain
\[
M_f(\beta+2r)=\sum_{s=0}^r\binom{\beta+2r}{\beta+2s}(k/2-1)^{2(r-s)}C_f(\beta+2s),
\qquad 0\le r<q.
\]
Thus the two moment maps are related by a unit lower triangular matrix and have the same rank, proving the asserted bounds.
\end{proof}

\subsection{Arithmetic consequences}\label{subsection_period-arithmetic}

To conclude this section, we study rational spans of the periods $r_n(f)=r_n(f,0)$ of a single newform. Shimura's period relations identify suitably normalized periods of its Galois conjugates with the conjugates of algebraic numbers. This allows the preceding interpolation results to control the rational span of periods of one parity.

Throughout this subsection, let $f$ be a normalized newform of weight $k$ and level $N\ge1$, with trivial nebentypus and standard Fourier expansion $f(z)=\sum_{n\ge1}a_f(n)e^{2\pi inz}$. Write
\[
K_f=\Q(a_f(n):n\ge1),\qquad e=[K_f:\Q].
\]
Choose $\Gamma$ and $\chi$ from Proposition~\ref{proposition_groups} for this same $N$ so that $S_k(\Gamma,\chi)$ contains all Galois conjugates of $f$, and write
\[
d=d_k(\Gamma,\chi),\qquad \delta=g(\Gamma)-\ind_{\{k=2,\ \chi=1\}}.
\]
Fix distinct indices $n_1,\ldots,n_q\in\{0,\ldots,k-2\}$ of the same parity. If the chosen group $\Gamma$ contains $W_N$, require $n_i+n_j\ne k-2$ for $i\ne j$, and exclude $n_j+1=k/2$ when $i^k\chi(W_N)=-1$.

By Shimura's period relations \cite[Theorem~1]{shimura1977periods}, choose a nonzero period $\Omega_f$ for this parity and set
\[
\xi_j=\frac{r_{n_j}(f)}{\Omega_f}\in K_f,\qquad 1\le j\le q.
\]

\begin{corollary}\label{corollary_arithmetic}
We have
\begin{equation}\label{equation_rational-span}
\min\{e,[e+q-d-\delta]_+\}
\le\dim_{\Q}\Span_{\Q}\{\xi_1,\ldots,\xi_q\}
\le\min\{e,q\}.
\end{equation}
In particular, a sufficient condition for the $\xi_j$ to span $K_f$ over $\Q$ is $q\ge d+\delta$. If $\delta=0$ and $e=d$, then the $\xi_j$, and hence the periods $r_{n_j}(f)$, are linearly independent over $\Q$ whenever $q\le e$.
\end{corollary}

\begin{proof}
Let $\sigma_1,\ldots,\sigma_e$ be the distinct embeddings of $K_f$ into $\C$, and set
\[
 W_f=\Span_\C\{f^{\sigma_1},\ldots,f^{\sigma_e}\}.
\]
The Hecke eigenvalues generate $K_f$, so the conjugates have distinct joint Hecke eigencharacters and form a basis of this $e$-dimensional subspace. Let $\mathcal P(F)=(r_{n_1}(F),\ldots,r_{n_q}(F))$ on the chosen ambient space. By \eqref{equation_period-definition} and Theorem~\ref{theorem_modular}, its rank is at least $\min\{d,[q-\delta]_+\}$. Restriction to $W_f$ lowers rank by at most its codimension $d-e$, giving
\[
\rank_\C(\mathcal P|_{W_f})\ge[\rank_\C\mathcal P-(d-e)]_+\ge\min\{e,[e+q-d-\delta]_+\}.
\]

We next use Shimura's relations to interpret this complex rank. The identity
\[
r_n(f)=(-1)^{n+1}n!\,\frac{L(f,n+1)}{(2\pi i)^{n+1}}
\]
shows that our periods differ from the standard Shimura normalization only by nonzero rational factors. The Galois equivariance in Shimura's period relations, valid also in weight two, provides nonzero periods $\Omega_i$ for the conjugate forms such that
\begin{equation}\label{equation_galois-period-matrix}
\sigma_i(\xi_j)=\frac{r_{n_j}(f^{\sigma_i})}{\Omega_i} \qquad(1\le i\le e,\ 1\le j\le q).
\end{equation}
Because the indices have the same parity, the same period $\Omega_i$ can be used for every $j$. Since the conjugates form a basis of $W_f$, dividing the rows of the transposed period matrix by $\Omega_i$ identifies its rank with that of $[\sigma_i(\xi_j)]_{i,j}$.

Since $K_f/\Q$ is separable, the embeddings $\sigma_i$ induce an isomorphism $K_f\otimes_{\Q}\C\simeq\C^e$. Under this isomorphism, the complexification of $\Span_{\Q}\{\xi_1,\ldots,\xi_q\}$ maps onto the column space of $[\sigma_i(\xi_j)]_{i,j}$. Hence
\[
 \rank_\C(\mathcal P|_{W_f})
 =\rank_\C[\sigma_i(\xi_j)]_{i,j}
 =\dim_\Q\Span_\Q\{\xi_1,\ldots,\xi_q\}.
\]
Together with the preceding rank bound, this proves \eqref{equation_rational-span} and its consequences. Multiplication by $\Omega_f\ne0$ preserves rational linear relations.

For the Atkin--Lehner choice of ambient space, the character on $\Gamma_0(N)$ is trivial and $k$ is even. Under these assumptions, the operators $W_e$ preserve forms with rational Fourier coefficients \cite[Corollary~5.1(2)]{brunault2020fourier}. Their matrices in a basis with rational Fourier coefficients therefore have rational entries, so they commute with Galois conjugation of Fourier coefficients. Consequently, for $W_e\in A$,
\[
(f^\sigma)|_kW_e=(f|_kW_e)^\sigma=\chi(W_e)f^\sigma,
\]
where $\chi(W_e)\in\{\pm1\}$. Thus every Galois conjugate of $f$ belongs to the same Atkin--Lehner character space.
\end{proof}

The ratios of the selected periods belong to $K_f$ whenever their denominators are nonzero, since $r_{n_j}(f)/r_{n_{j_0}}(f)=\xi_j/\xi_{j_0}$. The preceding rank bound gives a sufficient condition for these ratios to generate $K_f$.

\begin{corollary}\label{corollary_period-ratios}
If $q>d+\delta-e/2$, then at least one of the selected periods is nonzero, and for any $j_0$ with $r_{n_{j_0}}(f)\ne0$,
\[
K_f=\Q\left(\frac{r_{n_j}(f)}{r_{n_{j_0}}(f)}:1\le j\le q\right).
\]
In particular, the conclusion holds when $\delta=0$, $e=d$, and $q>d/2$.
\end{corollary}

\begin{proof}
By Corollary~\ref{corollary_arithmetic},
\[
\dim_\Q\Span_\Q\{\xi_1,\ldots,\xi_q\}
\ge\min\{e,[e+q-d-\delta]_+\}>e/2.
\]
For any $\xi_{j_0}\ne0$, let $F=\Q(\xi_j/\xi_{j_0}:1\le j\le q)\subseteq K_f$. Division by $\xi_{j_0}$ preserves the rational rank, so $[F:\Q]>e/2$. Every proper subfield of $K_f$ has degree at most $e/2$ by the tower law. Thus $F=K_f$, and $\xi_j/\xi_{j_0}=r_{n_j}(f)/r_{n_{j_0}}(f)$.
\end{proof}

\section{Critical period maps}\label{section_periods}

In this section, we study finite families of period functionals defined by integration along paths joining cusps. These include $r_n(\,\cdot\,,\alpha)$, obtained by integration along the vertical path from $\alpha$ to $i\infty$, for $\alpha\in\Q$ and $0\le n\le k-2$. The following proposition shows that finitely many of these functionals span the dual of the cusp-form space.

\begin{proposition}\label{proposition_hecke-periods}
Let $N\ge1$ be an integer and let $k\ge4$ be even. Let
\[
\mu_N=[\SL_2(\Z):\Gamma_0(N)],\qquad
B=\left\lceil\frac{k\mu_N}{12}\right\rceil.
\]
For $1\le q\le B$, define
\[
\mathcal A_q(f)=q^{k-2}\sum_{a=0}^{q-1}r_{k-2}(f,a/q).
\]
Then $f\longmapsto(\mathcal A_1(f),\ldots,\mathcal A_B(f))$ is injective on $S_k(\Gamma_0(N))$. Consequently, the period functionals $r_{k-2}(\,\cdot\,,a/q)$, with $1\le q\le B$ and $0\le a<q$, span its dual.
\end{proposition}
\begin{proof}
Let $\mathbb T\subset\operatorname{End}_{\C}(S_k(\Gamma_0(N)))$ be the full commutative Hecke algebra, including $T_p=U_p$ for $p\mid N$. Consider the space
\[
K=\{f\in S_k(\Gamma_0(N)):r_{k-2}(Tf,0)=0\text{ for every }T\in\mathbb T\}.
\]
Commutativity makes $K$ invariant under $\mathbb T$. If $K\ne0$, it contains a nonzero simultaneous eigenvector $g$ for all $T_n$. The identities $a_g(n)=a_1(T_ng)$ imply $a_g(1)\ne0$, so we normalize $a_g(1)=1$. The Euler product of $L(g,s)$ is absolutely convergent and nonzero for $\operatorname{Re}s>(k+1)/2$. Since $k-1>(k+1)/2$, it follows that
\[
r_{k-2}(g,0)=\frac{i^{k-1}(k-2)!}{(2\pi)^{k-1}}L(g,k-1)\ne0,
\]
a contradiction. Hence $K=0$ and the functionals $f\mapsto r_{k-2}(Tf,0)$, $T\in\mathbb T$, span the dual.

The Sturm bound and the perfect pairing $(T,f)\mapsto a_1(Tf)$ imply that $T_1,\ldots,T_B$ span $\mathbb T$ as a vector space; see \cite[Theorem~9.23]{stein2007modular}. Thus the functionals $f\mapsto r_{k-2}(T_nf,0)$, $1\le n\le B$, already span the dual of $S_k(\Gamma_0(N))$. Substituting the standard Hecke formula
\[
T_nf(z)=\sum_{\substack{uv=n\\(u,N)=1}}\frac{u^{k-1}}{v}
             \sum_{a=0}^{v-1}f\left(\frac{uz+a}{v}\right)
\]
into the period integral and changing variables, we obtain
\[
r_{k-2}(T_nf,0)=\sum_{\substack{uv=n\\(u,N)=1}}\mathcal A_v(f).
\]
For $1\le n\le B$, these identities define a triangular change of coordinates with diagonal entries $1$. This change preserves rank, so the functionals $\mathcal A_1,\ldots,\mathcal A_B$ also span the dual.
\end{proof}

The functionals $r_{k-2}(\,\cdot\,,a/q)$ supplied by Proposition~\ref{proposition_hecke-periods} involve an increasing number of rational starting points as the weight grows. We now use Chebyshev systems to obtain smaller spanning families and, for certain cusp-form spaces, explicit dual bases using paths independent of the weight. The following example illustrates this improvement.

\begin{example}\label{remark_period-level27}
At level $27$ and in even weight $k\ge4$, Proposition~\ref{proposition_hecke-periods} uses $B=3k$ functionals $\mathcal A_q$, formed from periods with rational starting points of denominator at most $3k$. In Section~\ref{subsection_cyclic-periods}, we select the following $3k-6$ period functionals with $\alpha=0,1/3,2/3$ as a basis of $S_k(\Gamma_0(27))^*$:
\begin{equation}\label{equation_period-level27-basis}
\begin{aligned}
&\{r_n(\,\cdot\,,0):0\le n\le k-2\}\\
&\cup\{r_n(\,\cdot\,,1/3):0\le n\le k-3\}\cup\{r_n(\,\cdot\,,2/3):0\le n\le k-4\}.
\end{aligned}
\end{equation}
For $k=2$, the single functional $r_0(\,\cdot\,,0)$ is a basis. For $k\ge6$, at least three values of $\alpha$ are needed for the functionals $r_n(\,\cdot\,,\alpha)$ to span the dual.
\end{example}

To study these period functionals, we count sign changes along real arcs. A form with real Fourier coefficients takes real values on the vertical paths through $0$ and $1/2$. After removing powers of $i$, its integrals on these paths are real moments against powers of $y$, and vanishing moments force sign changes. Suppose the arcs project injectively to pairwise disjoint subsets of the modular curve and avoid elliptic points. Lemma~\ref{lemma_divisor} then bounds the total number of sign changes in terms of the dimension and genus. Proposition~\ref{proposition_joint-rank} converts this bound into rank estimates over $\R$, which remain valid over $\C$ by Lemma~\ref{lemma_complexification}. The same argument applies to other parametrized paths with a common real form.

We first apply this principle directly to the standard real locus of $X_0(N)$. Cyclic translation (Section \ref{subsection_cyclic-periods}) and decomposition into Fricke eigenspaces (Section \ref{subsection_AL-periods}) then allow us to work on smaller full cusp-form spaces and quotient curves, sharpening the bounds and yielding explicit bases in genus zero. We finish by comparing the scope of these constructions with that of the direct period families at $0$ and at $0,1/2$ in Table~\ref{table_small-levels}.

\subsection{Periods on the standard real locus}\label{subsection_real-periods}

Let $N\ge4$ and let $k\ge2$ be even. The standard real structure is induced by $\mathfrak c(z)=-\overline z$. Let $b_N$ count the open arcs obtained by removing the real
cusps from its real locus. Each real component is a circle containing a cusp, so $b_N$ equals the number of real cusps. For $N=2^aM$ with $M$ odd, the cusp count in \cite[Corollary~5.1.4 and Proposition~6.5.1]{snowden2011real} is
\begin{equation}\label{equation_real-arc-count}
b_N=2^{\omega(M)}
\begin{cases}
1,&a=0,\\
2,&a=1,\\
3,&a=2,\\
4,&a\ge3.
\end{cases}
\end{equation}
Here $\omega(M)$ counts the distinct prime divisors of $M$, with $\omega(1)=0$. For each arc, choose a geodesic lift $A_j(i\R_{>0})$, with $A_j\in\mathrm{GL}_2^+(\Q)$, and define the associated period integrals
\begin{equation}\label{equation_real-path-periods}
 p_{j,n}(f)=\int_0^{i\infty}(f|_kA_j)(z)z^n\,dz,
 \qquad 1\le j\le b_N,\quad 0\le n\le k-2.
\end{equation}

\begin{proposition}\label{proposition_real-periods}
Let $\mathcal P$ be the map whose coordinates are $p_{j,n}$ for any $q$ distinct index pairs $(j,n)$ in \eqref{equation_real-path-periods}. Then
\begin{equation}\label{equation_real-period-rank}
\min\{d_k(N),[q-g_N+\ind_{\{k=2\}}]_+\}\le\rank\mathcal P\le\min(d_k(N),q).
\end{equation}
In particular, a sufficient condition for the selected functionals to span the dual is $q\ge d_k(N)+g_N-\ind_{\{k=2\}}$. If $g_N=0$, their rank is exactly $\min(d_k(N),q)$. Each $p_{j,n}$ is a rational linear combination of the period functionals $r_m(\,\cdot\,,\alpha)$, $0\le m\le k-2$, where $\alpha\in\{A_j0,A_j\infty\}\cap\Q$.
\end{proposition}

\begin{proof}
By \cite[Proposition~5.3.1]{snowden2011real}, the standard real locus contains no elliptic points of order two. To exclude order three, suppose that $\tau\xi$ represents
such a real elliptic point, where $\tau\in\SL_2(\Z)$ and $\xi=(-1+i\sqrt3)/2$. Set
\[
\rho=\begin{pmatrix}0&-1\\1&1\end{pmatrix},\qquad T=\begin{pmatrix}1&1\\0&1\end{pmatrix}.
\]
The stabilizer of $\tau\xi$ in $\SL_2(\Z)$ is generated by $\tau\rho\tau^{-1}$ and is contained in $\Gamma_0(N)$. Since $\mathfrak c(\tau\xi)=\tau^\#T\xi$, reality then implies $\tau^\#T\tau^{-1}\in\Gamma_0(N)$. Their lower-left entries must be divisible by $N$, so
\[
N\mid u^2-uv+v^2,\qquad N\mid u(u-2v),
\]
where $(u,v)$ is the bottom row of $\tau$. Since $\gcd(u,v)=1$, the first divisibility condition makes both $u$ and $v$ units modulo $N$. The second condition forces $u\equiv2v\pmod N$, and substitution into the first implies $N\mid3$, contrary to $N\ge4$. Thus the standard real locus contains no elliptic points. With elliptic points excluded, Snowden's description \cite[Sections~3.3--3.5 and Proposition~5.1.3]{snowden2011real} shows that these open arcs lift to geodesics with endpoints in $\Q\cup\{\infty\}$ and interiors that project injectively onto the corresponding pairwise disjoint arcs.

The pointwise reflection in the lifted real arc $A_j(i\R_{>0})$ is
\[
A_j\mathfrak c A_j^{-1}
=A_j(A_j^\#)^{-1}\mathfrak c
\in\Gamma_0(N)\mathfrak c.
\]
It follows that $A_j^\#A_j^{-1}\in\Gamma_0(N)$. Consequently, for $f^c=f$, we have
\[
(f|_kA_j)^c=f|_kA_j^\#=f|_kA_j.
\]
Thus all $F_j(y)=(f|_kA_j)(iy)$ are real-valued on the same real form, namely the space of forms with real Fourier coefficients. For $f\ne0$, each $F_j$ is nonzero. Applying Lemma~\ref{lemma_divisor} to any finite collection of zeros on these disjoint arcs, we obtain
\[
\sum_j Z_{(0,\infty)}(F_j)\le d_k(N)+g_N-1-\ind_{\{k=2\}}.
\]
In particular, the zero sets are finite. By cuspidality at $A_j0$ and $A_j\infty$, each $F_j$ decays exponentially at infinity and satisfies $F_j(y)=O(y^{-k}e^{-c/y})$ at zero. All auxiliary real-power moments therefore converge absolutely. Any increasing sequence of distinct powers of $y$ is an ECT-system on $(0,\infty)$. Since
\[
i^{-n-1}p_{j,n}(f)=\int_0^\infty F_j(y)y^n\,dy,
\]
Proposition~\ref{proposition_joint-rank} with $\delta=g_N-\ind_{\{k=2\}}$ and Lemma~\ref{lemma_complexification} yield \eqref{equation_real-period-rank}. 

To express these integrals in terms of $r_m(f,\alpha)$, write $A=\left(\begin{smallmatrix}a&b\\c&e\end{smallmatrix}\right)$ and substitute $w=Az$ in the integral:
\begin{equation}\label{equation_rational-path-transfer}
\int_0^{i\infty}(f|_kA)(z)z^n\,dz=(\det A)^{1-k/2}\int_{A0}^{A\infty}f(w)(ew-b)^n(a-cw)^{k-2-n}\,dw.
\end{equation}
The polynomial has degree at most $k-2$ and rational coefficients. Splitting the path at $\infty$ and expanding this polynomial about each finite endpoint $\alpha$ expresses the integral as a rational linear combination of $r_m(f,\alpha)$, $0\le m\le k-2$. 
\end{proof}

The span of any selected family of $p_{j,n}$ is contained in the span of the period functionals $r_m(\,\cdot\,,\alpha)$, $0\le m\le k-2$, where $\alpha$ ranges over the finite endpoints of the corresponding lifts. Hence the latter family has at least the same rank and spans $S_k(\Gamma_0(N))^*$ whenever the selected $p_{j,n}$ do.

For the complete family in \eqref{equation_real-path-periods}, the number of coordinates is $q=b_N(k-1)$. Combining \eqref{equation_real-period-rank} with the dimension formula
\[
d_k(N)=\frac{\mu_N}{12}k+O_N(1),\qquad \mu_N=[\SL_2(\Z):\Gamma_0(N)],
\]
we obtain, for fixed $N$ as $k\to\infty$ through even integers,
\begin{equation}\label{equation_real-period-asymptotic}
 \rank\mathcal P
 =\min\{\mu_N/12,b_N\}\,k+O_N(1).
\end{equation}

\begin{example}\label{remark_period-level40}
At level $40$, we have $b_{40}=8$ and $g_{40}=3$. For even $k\ge4$, the dimension is $d_k(40)=6k-10$, and
\[
 b_{40}(k-1)-g_{40}=8k-11\ge d_k(40).
\]
For $k=2$, the rank lower bound is $\min\{3,8-3+1\}=3$. Thus the complete family $p_{j,n}$ on the eight standard real arcs spans $S_k(\Gamma_0(40))^*$ in every even weight. Consequently, the functionals $r_n(\,\cdot\,,\alpha)$, $0\le n\le k-2$, also span this dual when $\alpha$ ranges over the finite endpoints of chosen lifts of these eight arcs.
\end{example}

\subsection{Cyclic translation}\label{subsection_cyclic-periods}

We next use cyclic translation to decompose the period map at equally spaced rational points into maps on character spaces of lower level. Assume in this subsection that
\[
\ell\in\{2,3,4,6,8,12,24\},\qquad \ell^2\mid N,\qquad M=N/\ell\ge4,\qquad k\ge2\ \text{even}.
\]
Let $\zeta=e^{2\pi i/\ell}$. For $0\le r<\ell$, set
\[
\chi_r\begin{pmatrix}a&b\\c&d\end{pmatrix}=\zeta^{rab},
\qquad d_r=d_k(\Gamma_0(M),\chi_r).
\]
The $\chi_r$ are characters of the group $\Gamma_0(M)$. Write $\mathcal R=\{r:d_r>0\}$. Maxima over empty index sets are taken to be zero.

The following theorem bounds the ranks in terms of the dimensions $d_r$ and the genus $g_M$. If the genus is zero, then the bounds coincide.

\begin{theorem}\label{theorem_cyclic-periods}
Let $I,J\subset\{0,\ldots,k-2\}$ and $q=|I|+|J|$. Let $\mathcal P$ be the map on $S_k(\Gamma_0(N))$ whose coordinates are
\[
\{r_n(\,\cdot\,,j/\ell):n\in I,\ 0\le j<\ell\}
\ \cup\
\{r_n(\,\cdot\,,(2j+1)/(2\ell)):n\in J,\ 0\le j<\ell\}.
\]
Then
\[
\sum_{r\in\mathcal R}\min\{d_r,[q-g_M+\ind_{\{k=2,r=0\}}]_+\}
\le\rank\mathcal P
\le\sum_{r\in\mathcal R}\min(d_r,q).
\]
A necessary condition for spanning is $q\ge\max_r d_r$, and a sufficient condition is $q\ge\max_{r\in\mathcal R}(d_r+g_M-\ind_{\{k=2,r=0\}})$. If $g_M=0$, the rank is exactly $\sum_r\min(d_r,q)$, and the selected functionals span the dual if and only if $q\ge\max_r d_r$.
\end{theorem}

The theorem treats complete translation orbits at the chosen indices. The next proposition selects individual periods from these orbits, giving a basis in genus zero.

\begin{proposition}\label{proposition_cyclic-selection}
Let $t_r=d_r+g_M-\ind_{\{k=2,r=0\}}$ for $r\in\mathcal R$, $t_{\max}=\max_{r\in\mathcal R}t_r$, and $s_j=\#\{r\in\mathcal R:t_r\ge j\}$. Choose either $(n_j,\alpha_j)=(j-1,0)$ or the alternating enumeration
\[
(n_{2a+1},\alpha_{2a+1})=(a,0),\qquad
(n_{2a+2},\alpha_{2a+2})=(a,1/(2\ell)),\qquad a\ge0.
\]
Assume $t_{\max}\le k-1$ in the first enumeration, or $t_{\max}\le2(k-1)$ in the second. Then the period functionals
\[
\bigcup_{j=1}^{t_{\max}}
\{r_{n_j}(\,\cdot\,,\alpha_j+v/\ell):0\le v<s_j\}
\]
span $S_k(\Gamma_0(N))^*$. Their number is $\sum_{r\in\mathcal R}t_r$. If $g_M=0$, this family is a basis.
\end{proposition}

For the complete critical families, each component has $q=k-1$ moments when one orbit of rational parameters is used, and $q=2(k-1)$ when two orbits are used. A sufficient condition for spanning is $t_{\max}\le k-1$ or $t_{\max}\le2(k-1)$, respectively. In genus zero, the family spans if and only if the corresponding inequality holds.

The proofs use the following geometric properties of the two paths.

\begin{lemma}\label{lemma_two-period-paths}
Define $\ell_a=\{a+iy:y>0\}$ for $a=0,1/2$. For $M\ge4$, the paths $\ell_0$ and $\ell_{1/2}$ project injectively to $X_0(M)$, with disjoint images containing no elliptic points.
\end{lemma}

\begin{proof}
The injectivity of $\ell_0$ and the absence of elliptic points on its image follow from the proof of Proposition~\ref{proposition_groups}. We first establish these properties for $\ell_{1/2}$.

Suppose that $\gamma=\left(\begin{smallmatrix}a&b\\c&d\end{smallmatrix}\right)
\in\Gamma_0(M)$ satisfies
\[
\gamma(1/2+iy)=1/2+iy', \qquad y,y'>0.
\]
If $c=0$, then $\gamma$ acts by an integral translation, and the equality of real parts forces $\gamma=\pm I$. Assume henceforth that $c\ne0$, so that $|c|\ge M\ge4$. Conjugating by the translation $T(z)=z+1/2$, write
\[
T^{-1}\gamma T=
\begin{pmatrix}a'&b'\\c&d'\end{pmatrix},
\qquad
a'=a-\frac c2,\quad
b'=b+\frac{a-d}{2}-\frac c4,\quad
d'=d+\frac c2.
\]
Since $T^{-1}\gamma T$ sends $iy$ to $iy'$, we have
\[
a'cy^2+b'd'=0, \qquad a'd'-b'c=1.
\]
Eliminating $a'$ gives
\[
b'c=-\frac{c^2y^2}{c^2y^2+(d')^2}\in[-1,0).
\]
Since $4b'\in\mathbb Z$, the nonzero quantity $b'c$ has absolute value at least $|c|/4\ge1$. Thus $b'c=-1$, which forces $d'=0$, and then $a'=0$. Consequently $c=2a$ and $d=-a$, so
\[
1=ad-bc=-a(a+2b).
\]
Hence $|a|=1$ and $|c|=2$, a contradiction. Therefore only $\pm I$ can identify points on $\ell_{1/2}$. This proves injectivity and also excludes nontrivial elliptic stabilizers on this path.

It remains to prove that the two images are disjoint. Suppose that $\gamma\in\Gamma_0(M)$ satisfies
\[
\gamma(iy)=1/2+iy',
\qquad y,y'>0.
\]
An integral translation cannot give such an identification, so $c\ne0$. Also $d\ne0$, since $d=0$ would imply $bc=-1$, contrary to $M\mid c$.
Set $A=2a-c$ and $B=2b-d$. The real-part condition and $ad-bc=1$ give
\[
cAy^2+dB=0, \qquad Ad-Bc=2.
\]
Eliminating $A$, we obtain
\[
Bc=-\frac{2c^2y^2}{c^2y^2+d^2}\in(-2,0).
\]
As $Bc$ is an integer, it must equal $-1$. This contradicts $M\mid c$ and completes the proof.
\end{proof}

For the character spaces used here, cuspidality at $0$, $1/2$, and $\infty$ ensures absolute convergence of the moments against every real power of $y$, including the auxiliary moments in Proposition~\ref{proposition_joint-rank}. On any nonzero full space $S_k(\Gamma,\chi)$, Lemma~\ref{lemma_divisor} bounds the total number of zeros on arcs projecting injectively to pairwise disjoint subsets of $X_\Gamma$ and avoiding elliptic points by $d_k(\Gamma,\chi)+\delta-1$, where
\begin{equation}\label{equation_joint-genus}
\delta=g(\Gamma)-\ind_{\{k=2,\ \chi=1\}}\ge0.
\end{equation}
We will apply this bound after verifying a common real form on the relevant arcs.

\begin{proof}[Proof of Theorem~\ref{theorem_cyclic-periods}]
For $f\in S_k(\Gamma_0(N))$, let
\begin{equation}\label{equation_cyclic-projector}
 f_r(z)=\frac1\ell\sum_{j=0}^{\ell-1}\zeta^{-rj}f(z+j/\ell)
       =\sum_{\substack{n\ge1\\ n\equiv r\pmod\ell}}
          a_f(n)e^{2\pi inz}.
\end{equation}
For the allowed values of $\ell$, every integer $a$ coprime to $\ell$ satisfies $a^2\equiv1\pmod\ell$. Since conjugation changes the upper-right entry by $(a-d)/\ell-c/\ell^2$, an integer, translation by $1/\ell$ normalizes $\Gamma_0(N)$. Thus each $f_r$ belongs to $S_k(\Gamma_0(N))$. On $\Gamma_0(M)$, the map $\psi(\gamma)=ab\bmod\ell$ is a surjective
homomorphism, with
\[
\ker\psi=\{\gamma\in\Gamma_0(M):\ell\mid b\}=H\Gamma_0(N)H^{-1},\qquad H=\operatorname{diag}(\sqrt\ell,1/\sqrt\ell).
\]
The homomorphism property follows by matrix multiplication using $\ell\mid M$ and the congruence $a^2\equiv1\pmod\ell$, and surjectivity follows from $\psi\left(\begin{smallmatrix}1&1\\0&1\end{smallmatrix}\right)=1$.
Since $f_r(z+1/\ell)=\zeta^rf_r(z)$, the scaled form $z\mapsto f_r(z/\ell)$ belongs to $S_k(\Gamma_0(M),\chi_r)$. Thus scaling induces the isomorphism
\[
S_k(\Gamma_0(N))\longrightarrow
\bigoplus_{r=0}^{\ell-1}S_k(\Gamma_0(M),\chi_r),\qquad
f\longmapsto\bigl(f_r(\,\cdot\,/\ell)\bigr)_r.
\]
For surjectivity, note that every $F\in S_k(\Gamma_0(M),\chi_r)$ is obtained by scaling the $r$th Fourier summand $F(\ell z)$ at level $N$. Positive cusp exponents are multiplied by ramification indices, so cuspidality is preserved in both directions. In particular, $\sum_rd_r=d_k(N)$.

Since $\chi_r(\gamma^\#)=\overline{\chi_r(\gamma)}$, coefficient conjugation preserves each space, even for nonreal $\chi_r$. By \eqref{equation_cyclic-projector}, conjugation also commutes with the Fourier projections: $(f_r)^c=(f^c)_r$. For $F=F^c\in S_k(\Gamma_0(M),\chi_r)$, the expansion
\[
F(z)=\sum_{j+r/\ell>0}a_j e^{2\pi i(j+r/\ell)z},\qquad a_j\in\R,
\]
implies $F(iy)\in\R$ and $e^{-\pi ir/\ell}F(1/2+iy)\in\R$. Thus the period functionals on the two paths become real-valued on this same real form after multiplication by $i^{-n-1}$ and $i^{-n-1}e^{-\pi ir/\ell}$, respectively. By Lemma~\ref{lemma_two-period-paths}, the two paths project injectively to disjoint subsets of $X_0(M)$ containing no elliptic points. Applying Lemma~\ref{lemma_divisor} to the zeros of any nonzero $F$ on both paths therefore bounds their total number, counted with multiplicity, by $d_r+g_M-\ind_{\{k=2,r=0\}}-1$. On each path, powers of $y$ with distinct real exponents form an ECT-system on $(0,\infty)$ by Lemma~\ref{lemma_kernels}(1). Cuspidality at both endpoints ensures absolute convergence of the moments for every real exponent. We may therefore adjoin moments with further distinct exponents on either path, preserving the ECT property and absolute convergence, as required
by Proposition~\ref{proposition_joint-rank}.

For $\alpha=0$ or $1/(2\ell)$, translation and dilation give
\begin{equation}\label{equation_period-DFT}
\ell^{-1}\sum_{j=0}^{\ell-1}\zeta^{-rj}r_n(f,\alpha+j/\ell)=r_n(f_r,\alpha)
=\ell^{-n-1}r_n\bigl(f_r(\,\cdot\,/\ell),\ell\alpha\bigr).
\end{equation}
By \eqref{equation_period-DFT}, the Fourier-transformed rows indexed
by $r$ depend only on $f_r(\,\cdot\,/\ell)$. Up to nonzero row
factors, they are the $q=|I|+|J|$ moments
\[
\bigl(r_n(f_r(\,\cdot\,/\ell),0)\bigr)_{n\in I},\qquad
\bigl(r_n(f_r(\,\cdot\,/\ell),1/2)\bigr)_{n\in J}.
\]
Thus, under the isomorphism $f\mapsto(f_r(\,\cdot\,/\ell))_r$, the transformed period map is a direct sum of moment maps on $S_k(\Gamma_0(M),\chi_r)$.
The summands with $d_r=0$ contribute only zero rows. After omitting these rows and applying the phase normalizations above, the blocks are complexifications of real moment maps. Apply \eqref{equation_joint-rank} to their direct sum over $\R$, with $q_r=q$ and
$\delta_r=g_M-\ind_{\{k=2,r=0\}}$ for $r\in\mathcal R$.
Complexification preserves rank by Lemma~\ref{lemma_complexification}, as do the invertible row operations, so these are also the asserted bounds for $\mathcal P$.
If $g_M=0$, the exceptional case $k=2$, $r=0$ has
$d_0=g_M=0$. Hence $\delta_r=0$ on every nonzero summand,
and the rank bounds coincide.
\end{proof}

\begin{proof}[Proof of Proposition~\ref{proposition_cyclic-selection}]
Suppose the selected periods vanish. Inductively, assume $f_r=0$ for $t_r<j$. Since $f=\sum_rf_r$, the next equations are
\[
 0=r_{n_j}(f,\alpha_j+v/\ell)
   =\sum_{\substack{r\in\mathcal R\\t_r\ge j}}
       \zeta^{rv}r_{n_j}(f_r,\alpha_j),\qquad 0\le v<s_j.
\]
Their coefficient matrix is Vandermonde with distinct nodes $\zeta^r$. Hence $r_{n_j}(f_r,\alpha_j)=0$ for every $r$ with $t_r\ge j$. At $j=t_r$, these $t_r$ moment equations force $f_r=0$ by the injectivity of the complexified component maps established in the proof of Theorem~\ref{theorem_cyclic-periods}, together with the dilation in \eqref{equation_period-DFT}. This proves injectivity.

In either enumeration, the selected pairs $(n,\alpha)$ are distinct, and the bounds on $t_{\max}$ place all indices in the critical range. The number of selected periods is $\sum_{j=1}^{t_{\max}}s_j=\sum_{r\in\mathcal R}t_r$. When $g_M=0$, we have $t_r=d_r$ on $\mathcal R$. The injective map then has exactly $d_k(N)$ coordinates, so the selected periods form a dual basis.
\end{proof}

\begin{proof}[Verification of Example~\ref{remark_period-level27}]
Take $\ell=3$ and $M=9$. The curve $X_0(9)$ has genus zero, index $12$, four cusps, and no elliptic points. For a reduced cusp $a/c$, including $1/0=\infty$, its width and a positively oriented primitive stabilizer are
\[
w=\frac{9}{\gcd(c^2,9)},\qquad\gamma_{a/c}=
\begin{pmatrix}1-acw&a^2w\\-c^2w&1+acw\end{pmatrix}.
\]
Let $\alpha_{r,s}\in(0,1]$ be determined by $\chi_r(\gamma_s)=e^{2\pi i\alpha_{r,s}}$ at each cusp $s$, and put $b_r=\sum_s\alpha_{r,s}$. In the local parameter $q_s$, a cusp form has expansion $q_s^{\alpha_{r,s}}h(q_s)$ with $h$ holomorphic. Its order as a section of the cusp-form line bundle is therefore $\ord_s f-\alpha_{r,s}$. For a nonzero meromorphic section $f$, the form $f^3$ has trivial character and weight $3k$. Applying the valence formula to $f^3$ and dividing the resulting identity by three, we find
\[
\sum_{P\in\Gamma_0(9)\backslash\HH}\ord_P f+\sum_s\ord_s f=k,
\]
since the index is $12$ and there are no elliptic points. Thus the line bundle has degree $k-b_r$. We apply the Riemann--Roch theorem on $X_0(9)\simeq\mathbb P^1$ and then get $d_r=[k+1-b_r]_+$. At $0,1/3,2/3,\infty$, the stabilizers above have $\chi_1$-values $1,\zeta,\zeta,\zeta$. Hence
\[
(b_0,b_1,b_2)=(4,2,3),\qquad(d_0,d_1,d_2)=([k-3]_+,k-1,k-2).
\]
With these dimensions, the spanning criterion in Theorem~\ref{theorem_cyclic-periods} is satisfied. The first enumeration in Proposition~\ref{proposition_cyclic-selection} selects precisely the basis in \eqref{equation_period-level27-basis} and, in weight two, the single functional stated in the example. For $k\ge6$, at most $2(k-1)$ critical periods are available at two parameters, whereas $d_k(27)=3k-6>2(k-1)$. This proves minimality.
\end{proof}

\subsection{Atkin--Lehner translates}\label{subsection_AL-periods}

Assume that $N\ge5$ is odd and $k=2m\ge2$. Atkin--Lehner translates give explicit representatives of the standard real arcs considered in Section~\ref{subsection_real-periods}. We use the Fricke decomposition to estimate the rank on each sign component in terms of the genus of the Fricke quotient. By \eqref{equation_real-arc-count}, $b_N=2^{\omega(N)}$. Let $\Gamma^+=\Gamma_0(N)\cup\Gamma_0(N)W_N$, so $g(\Gamma^+)=h_N$. For $\tau=\pm1$, let $\chi_\tau$ be trivial on $\Gamma_0(N)$ with $\chi_\tau(W_N)=(-1)^m\tau$. Index the Fricke eigenspaces by the sign of the functional equation:
\[
S_k(\Gamma_0(N))=\bigoplus_{\tau=\pm1}S_k(\Gamma^+,\chi_\tau),\qquad d_\tau=d_k(\Gamma^+,\chi_\tau).
\]
Thus $d_+$ and $d_-$ are the dimensions of the components with functional-equation signs $+1$ and $-1$, respectively. In the notation of \eqref{equation_parity}, $\tau=(-1)^\beta$, while the Fricke eigenvalue is $(-1)^m\tau$.

Choose one exact divisor $e\parallel N$ (that is, $e\mid N$ and $\gcd(e,N/e)=1$) from each pair $\{e,N/e\}$, taking $1$ from $\{1,N\}$. There are $b_N/2$ such pairs. For a chosen $1<e<N$, choose $1\le t_e<N/e$ with $et_e\equiv1\pmod{N/e}$, and set
\begin{equation}\label{equation_period-AL-matrix}
b_e=e(et_e-1)/N,\qquad
W_e=\frac1{\sqrt e}\begin{pmatrix}e&b_e\\N&et_e\end{pmatrix},
\end{equation}
and $W_1=I$. These are Atkin--Lehner representatives as in \eqref{equation_atkin-matrix}.

To each chosen $e$ and $a\in\{0,1/2\}$, we associate the critical period functionals on $W_e\ell_a$:
\[
f\longmapsto r_n(f|_kW_e,a),\qquad 0\le n\le k-2.
\]
Fix any $b$ of these $b_N$ paths, with $1\le b\le b_N$. The following theorem bounds the rank of the associated period map.

\begin{theorem}\label{theorem_AL-periods}
Let $\mathcal P$ be the map on $S_k(\Gamma_0(N))$ whose coordinates
are $r_n(\,\cdot\,|_kW_e,a)$, $0\le n\le k-2$, on the $b$ selected
paths. Let
\[
q_+=bm,\qquad q_-=b(m-1).
\]
Then
\begin{equation}\label{equation_AL-period-rank}
\sum_{\tau=\pm1}
\min\{d_\tau,[q_\tau-h_N+\ind_{\{k=2,\tau=-1\}}]_+\}
\le\rank\mathcal P
\le\sum_{\tau=\pm1}\min(d_\tau,q_\tau).
\end{equation}
A necessary condition for spanning is $q_\tau\ge d_\tau$ for $\tau=\pm1$, and a sufficient condition is $q_\tau\ge d_\tau+h_N-\ind_{\{k=2,\tau=-1\}}$ for every $\tau$ with $d_\tau>0$. If $h_N=0$, the upper and lower bounds coincide, and the complete family spans if and only if $d_+\le bm$ and $d_-\le b(m-1)$.
\end{theorem}
The same rank bounds and spanning criteria apply to subfamilies whose selected indices on each path are closed under $n\mapsto k-2-n$. Here $q_-$ counts the selected pairs with indices $n$ and $k-2-n$, where $0\le n<m-1$, and $q_+$ is this
number plus the number of selected periods with index $n=m-1$, which we call central periods.

In genus zero, whenever the spanning criterion is satisfied, the following selection extracts a basis directly from the periods $r_n(\,\cdot\,|_kW_e,a)$ on the chosen paths.

\begin{proposition}\label{proposition_AL-selection}
Suppose $h_N=0$, $d_+\le bm$, and $d_-\le b(m-1)$. Let $c=\min\{b,d_+\}$ and choose $c$ central periods on distinct paths. Order arbitrarily the $b(m-1)$ pairs of periods with indices $n$ and $k-2-n$ on the same path, where $0\le n<m-1$.
Using this ordering, select
\begin{enumerate}
\item the period with index $n$ from each of the first
$\max\{d_+-c,d_-\}$ pairs;
\item the period with index $k-2-n$ from each of the first
$\min\{d_+-c,d_-\}$ pairs.
\end{enumerate}
Together with the $c$ central periods, these form a basis of
$S_k(\Gamma_0(N))^*$.
\end{proposition}
We next express the functionals $r_n(\,\cdot\,|_kW_e,a)$ in terms of $r_j(\,\cdot\,,\alpha)$, with $\alpha\in\{W_ea,W_e\infty\}\cap\Q$. For $e>1$, apply \eqref{equation_rational-path-transfer} with
\[
A=\sqrt e\,W_e\begin{pmatrix}1&a\\0&1\end{pmatrix}\in\mathrm{GL}_2^+(\Q),
\qquad a=0,1/2.
\]
Multiplication by the positive scalar $\sqrt e$ does not change the normalized slash action, so each $r_n(\,\cdot\,|_kW_e,a)$ is a rational linear combination of $r_j(\,\cdot\,,\alpha)$, $0\le j\le k-2$, with $\alpha$ a finite endpoint of the path. The two paths for each $e>1$ share the endpoint $e/N$ and therefore contribute at most three finite endpoints. Together with $0,1/2$ for $e=1$, the $b_N$ paths use at most $2+3(b_N/2-1)=3b_N/2-1$ rational points. Thus the basis elements in Proposition~\ref{proposition_AL-selection} are explicit rational linear combinations of these period functionals. In weight two, $d_-=h_N$, $d_+=g_N-h_N$, and $q_-=0$. Thus the complete family on the $b$ chosen paths in Theorem~\ref{theorem_AL-periods} spans if and only if $h_N=0$ and $b\ge g_N$.

When all $b_N$ paths are used, Lemma~\ref{lemma_AL-period-arcs} identifies their images with the standard real arcs. An invertible change of polynomial coordinates on each path therefore identifies the complete period map with that of Section~\ref{subsection_real-periods}. The two maps have the same rank in every weight, so \eqref{equation_real-period-asymptotic} also gives the asymptotic rank here. Thus a necessary condition for the complete family to span in an unbounded sequence of even weights is $\mu_N\le12b_N$.

The Fricke decomposition gives separate control of the two sign components. For $k\ge4$, the upper and lower bounds in \eqref{equation_AL-period-rank} differ by at most $2h_N$, compared with at most $g_N$ in Proposition~\ref{proposition_real-periods}. When $h_N=0$, the rank formula is exact, and Proposition~\ref{proposition_AL-selection} selects a basis whenever the family spans. For example, at level $35$, Proposition~\ref{proposition_real-periods} gives only $d_k(35)-1\le\rank\mathcal P\le d_k(35)$ for $k\ge4$, whereas Theorem~\ref{theorem_AL-periods} proves spanning. The following example gives an explicit basis selection.

\begin{example}\label{remark_period-level35}
Choose $e=1,7$ in \eqref{equation_period-AL-matrix}, with $t_7=3$,
$b_7=4$. For every even $k\ge2$,
\[
\Span_\C\{r_n(\,\cdot\,,\alpha):0\le n\le k-2,\
\alpha\in\{0,1/2,1/5,4/21,15/77\}\}=S_k(\Gamma_0(35))^*.
\]
For $k\ge4$, delete any one pair of noncentral periods, with indices
$n,k-2-n$ on the same path, from the complete critical family on these four paths. The remaining $4k-6$ functionals form a dual basis. For $k=2$, any three of the four central periods form a basis. By \eqref{equation_rational-path-transfer}, these bases are explicit rational linear combinations of the displayed period functionals.
\end{example}

We begin the proofs with the geometry of the paths.

\begin{lemma}\label{lemma_AL-period-arcs}
The $b_N$ paths $W_e\ell_a$, $a=0,1/2$, project injectively to $X_0(N)$, with pairwise disjoint images containing no elliptic points. On the Fricke quotient $X_{\Gamma^+}=X_0(N)/\langle W_N\rangle$, their open half-arcs
\begin{equation}\label{equation_period-AL-halves}
W_e\{iy:y>N^{-1/2}\},\qquad
W_e\{1/2+iy:y>(2\sqrt N)^{-1}\}
\end{equation}
have these same properties. Every $W_e$ preserves each character space and its real form $f^c=f$. On both base paths the restrictions of $f|_kW_e$ are real on this real form.
\end{lemma}
\begin{proof}
We first treat the character spaces and their real forms, then the paths on $X_0(N)$, and finally their images on the Fricke quotient.

The matrices $W_e$ normalize $\Gamma_0(N)$ and commute with $W_N$ modulo $\Gamma_0(N)$, so they preserve each $S_k(\Gamma^+,\chi_\tau)$. For $e>1$, we have
\[
W_e^\#W_e^{-1}=\begin{pmatrix}
2et_e-1&-2b_e\\
-2Nt_e&2et_e-1
\end{pmatrix}\in\Gamma_0(N).
\]
Thus $(f|_kW_e)^c=f^c|_kW_e$; the same identity is immediate for $W_1=I$. Each $W_e$ therefore preserves the real form. For $f=f^c$, the Fourier coefficients of $f|_kW_e$ are real, so its values at $iy$ and $1/2+iy$ are real.

Each $W_e$ normalizes $\Gamma_0(N)$, so the injectivity and absence of elliptic points established in Lemma~\ref{lemma_two-period-paths} also hold for the translated paths. It remains to prove pairwise disjointness. We first distinguish the paths up to $\Gamma_0(N)$-equivalence. Represent a cusp $u/v$ in lowest terms, including $\infty=1/0$, by the primitive integer vector $(u,v)^t$. The cusp denominator $\gcd(v,N)$ is invariant under $\Gamma_0(N)$, with value $N$ at $\infty$. For two endpoints represented by $(u,v)^t$ and $(u',v')^t$, the absolute determinant $|uv'-u'v|$ is also invariant: both vectors are multiplied by the same matrix of determinant one. For $e>1$, the paths $W_e\ell_0$ and $W_e\ell_{1/2}$ share the endpoint $e/N$, represented by the primitive vector $(1,N/e)^t$. Their other endpoints are represented by
\[
(b_e,et_e)^t,\qquad (e+2b_e,e(N/e+2t_e))^t.
\]
The identity $et_e-b_eN/e=1$ shows that the absolute determinants of these vectors with $(1,N/e)^t$ are one and two, respectively. The first vector is therefore primitive. A common divisor of the coordinates of the second divides two; since $e+2b_e$ is odd, the second vector is primitive as well. Moreover, $\gcd(t_e,N/e)=1$ and $N/e$ is odd, so
\[
\gcd(et_e,N)=\gcd(e(N/e+2t_e),N)=e.
\]
Thus both paths have unordered cusp-denominator pair $\{e,N/e\}$, while the absolute determinants distinguish the two choices of $a$. For $e=1$, the endpoint pairs are $\{0,\infty\}$ and $\{1/2,\infty\}$; their cusp-denominator pair is $\{1,N\}$, and the absolute determinants are again one and two. Since the pairs $\{e,N/e\}$ are distinct for the chosen exact divisors, all $b_N$ paths represent distinct unoriented $\Gamma_0(N)$-classes.

We now rule out intersections between these classes. Let $\mathfrak c(z)=-\bar z$.
The reflections in the base paths are $z\mapsto2a-\bar z$, $a=0,1/2$, and belong to $\Gamma_0(N)\mathfrak c$. Using $W_e^\#W_e^{-1}\in\Gamma_0(N)$, we may write
\[
W_e\mathfrak cW_e^{-1}
=\bigl(W_e^\#W_e^{-1}\bigr)^{-1}\mathfrak c
\in\Gamma_0(N)\mathfrak c.
\]
Since $W_e$ normalizes $\Gamma_0(N)$, the reflections in the translated paths also belong to this coset. The same holds for their $\Gamma_0(N)$-translates, because $\mathfrak c$ normalizes $\Gamma_0(N)$. If two path images intersected, suitable lifts would meet at a point $z\in\HH$. Their reflections would then have product in $\Gamma_0(N)$ fixing $z$. Since the paths contain no elliptic points, this product is projectively trivial. The two reflections, and hence their fixed geodesics, must therefore coincide, contradicting the distinct classes above. This proves pairwise disjointness on $X_0(N)$.

Finally, we determine the additional identifications on the Fricke quotient. The involutions
\[
J_0=W_N,\qquad
J_{1/2}=\frac1{\sqrt N}
\begin{pmatrix}
 N&-(N+1)/2\\
 2N&-N
\end{pmatrix}
\]
belong to $\Gamma_0(N)W_N$. For $J_{1/2}$, this follows from the oddness of $N$ and the identity
\[
J_{1/2}W_N^{-1}
=\begin{pmatrix}
(N+1)/2&1\\
N&2
\end{pmatrix}\in\Gamma_0(N).
\]
Their actions on the base paths are
\[
J_0(iy)=\frac{i}{Ny},\qquad
J_{1/2}\left(\frac12+iy\right)
=\frac12+\frac{i}{4Ny}.
\]
Thus they reverse the two paths, with fixed points at heights $N^{-1/2}$ and $(2\sqrt N)^{-1}$, respectively. For $L=W_e\ell_a$, the involution $R_L=W_eJ_aW_e^{-1}$ likewise reverses $L$. Since $W_e$ normalizes $\Gamma_0(N)$ and commutes with $W_N$ modulo this group, $R_L$ also belongs to the Fricke coset.

Suppose $z\in L$ and $\gamma z\in L'$ for $\gamma\in\Gamma_0(N)W_N$. Then $\gamma R_L\in\Gamma_0(N)$ sends $R_Lz$ to $\gamma z$. The disjointness and injectivity already proved on $X_0(N)$ give, respectively,
\[
L'=L,\qquad \gamma z=R_Lz.
\]
Hence the only additional identifications are the reversals within each path.
The open half-arcs in \eqref{equation_period-AL-halves} select one point from each two-point orbit of $R_L$ and omit the fixed point. They therefore project injectively to the Fricke quotient with pairwise disjoint images. Any new elliptic stabilizer would force a point to be fixed by $R_L$, which has been excluded. Thus their images contain no elliptic points.
\end{proof}

\begin{proof}[Proof of Theorem~\ref{theorem_AL-periods}]
Let $\kappa_0=N$ and $\kappa_{1/2}=4N$. For $f\in S_k(\Gamma^+,\chi_\tau)$, the reversals in Lemma~\ref{lemma_AL-period-arcs} give
\[
(f|_kW_e)\left(a+\frac{i}{\kappa_a y}\right)=\tau\kappa_a^m y^{2m}(f|_kW_e)(a+iy).
\]
Indeed, the reversing matrices on the base paths have automorphy factors $i\sqrt{\kappa_a}\,y$ and character values $(-1)^m\tau$. Substitute $y=e^x/\sqrt{\kappa_a}$ in the period integral and split the resulting integral at $x=0$. Pairing the two halves as in Lemma~\ref{lemma_mellin}, we obtain
\[
\begin{aligned}
&\kappa_a^{(n+1)/2}i^{-n-1}r_n(f|_kW_e,a)\\
&\qquad=2\int_0^\infty e^{mx}
(f|_kW_e)(a+ie^x/\sqrt{\kappa_a})
\begin{cases}
\cosh((n+1-m)x),&\tau=+1,\\
\sinh((n+1-m)x),&\tau=-1
\end{cases}
\,dx.
\end{aligned}
\]
Lemma~\ref{lemma_AL-period-arcs} supplies a common real form and disjoint halves, so \eqref{equation_joint-genus} bounds the total number of zeros by
$d_\tau+h_N-\ind_{\{k=2,\tau=-1\}}-1$ whenever $d_\tau>0$. Absolute convergence follows from cusp decay, as in Lemma~\ref{lemma_mellin}.

Normalize each period $r_n(\,\cdot\,|_kW_e,a)$ by the factor $\kappa_a^{(n+1)/2}i^{-n-1}$. Taking the half-sum and half-difference of the normalized functionals in each pair with indices $n,k-2-n$ is an invertible row operation. By the preceding integral formula, these vanish on the $\tau=-1$ and $\tau=+1$ summands, respectively. A central period vanishes on the $\tau=-1$ summand. Thus each noncentral pair contributes one moment to each component, while a central period contributes only to the $+$ component. For the complete family, the resulting counts are $q_+=bm$ and $q_-=b(m-1)$. On each arc, order the resulting nodes $u=(n+1-m)^2$ increasingly and divide each functional on the $\tau=-1$ component by its nonzero factor $n+1-m$. By Lemma~\ref{lemma_kernels}(3), the kernels $K_0(u,x)$ and $K_1(u,x)$ form the required ECT-systems, including the node $u=0$ on the $\tau=+1$ component. Further distinct positive nodes extend these ECT-systems, with absolute convergence again supplied by cusp decay. We can now apply \eqref{equation_joint-rank} to the direct sum of the real forms, taking $\delta_\tau=h_N-\ind_{\{k=2,\tau=-1\}}$ on the nonzero components. The resulting rank bounds over $\R$ remain valid over $\C$ by Lemma~\ref{lemma_complexification}. If $h_N=0$, then $\delta_\tau=0$ whenever $d_\tau>0$, since the exceptional weight-two component has dimension zero. Finally, in weight two $\chi_-=1$, and the $W_N$-invariant differentials descend to $X_{\Gamma^+}$. Hence $d_-=h_N$ and $d_+=g_N-h_N$, proving the stated weight-two criterion.
\end{proof}

\begin{proof}[Proof of Proposition~\ref{proposition_AL-selection}]
The spanning inequalities and the choice of $c$ ensure that both $d_+-c\le b(m-1)$ and $d_-\le b(m-1)$ hold. The prescribed selection is therefore possible and contains $d_++d_-$ periods. Suppose they vanish on a form $f$. We use the genus-zero injectivity for distinct moments proved above.

If $d_+-c\le d_-$, the central periods and the $d_+-c$ complete pairs impose $d_+$ vanishing moments on the $+$ component, forcing it to be zero. The $d_-$ selected lower-index periods then impose $d_-$ vanishing moments on the $-$ component, which must also be zero. If $d_+-c>d_-$, the $d_-$ complete pairs first force the $-$ component to vanish. The $c$ central periods and the $d_+-c$ selected lower-index periods then force the $+$ component to vanish. In each case the nodes on any given half-arc are distinct, since they come from distinct pairs, with the additional node $0$ for a central period. Thus $f=0$, and the selected periods form a basis. At weight two only the central periods are needed.
\end{proof}

\begin{proof}[Verification of Example~\ref{remark_period-level35}]
The curve $X_0(35)$ has index $48$, four cusps, no elliptic points, and genus $3$. For this level, the dimension formula and the full Fricke trace formula \cite[Section~3.2, Proposition~3.1]{martin2023root} read
\[
d_k(35)=4k-6\quad(k\ge4),\qquad
\operatorname{tr}(W_{35})=4(-1)^{k/2}+\ind_{\{k=2\}}.
\]
The sign of the functional equation is $(-1)^{k/2}$ times the Fricke eigenvalue. Thus $(d_+,d_-)=(2k-1,2k-5)$ for $k\ge4$, and $(d_+,d_-)=(3,0)$ for $k=2$. In weight two, $d_-=h_{35}=0$. The four paths have $q_+=2k$ and $q_-=2k-4$, so the complete family spans. For $k\ge4$, deleting one pair of noncentral periods with indices $n,k-2-n$ on the same path leaves $q_+=2k-1=d_+$ and $q_-=2k-5=d_-$. The exact rank formula therefore shows that the remaining functionals are linearly independent and hence form a basis. In weight two the same formula applies
to any three central periods.

The paths with $e=7$ have the common endpoint $1/5$ and their other endpoints are $4/21$ and $15/77$. For $e=1$, the finite endpoints are $0$ and $1/2$. Applying \eqref{equation_rational-path-transfer} to these four paths proves the assertions concerning the period functionals at the five rational points.
\end{proof}

\subsection{Conclusions}\label{subsection_small-levels}

Table~\ref{table_small-levels} compares four constructions for $1\le N<50$: the functionals $r_n(\,\cdot\,,\alpha)$ with $\alpha=0$ or $\alpha\in\{0,1/2\}$, integration along the standard real arcs, cyclic translation, and Atkin--Lehner translates. Throughout the table, $k\ge2$ is even. A complete family contains all indices $0\le n\le k-2$. The first column specifies the family to which each rank or spanning assertion applies; the rows for the standard real arcs record selected additional conclusions.

For the direct family at $0$, we use Theorem~\ref{theorem_basic}. For $N\ge4$, Lemma~\ref{lemma_two-period-paths} and Proposition~\ref{proposition_joint-rank}, followed by Lemma~\ref{lemma_complexification}, show that a sufficient condition for the complete family at $0,1/2$ to span the dual is
\[
2(k-1)\ge d_k(N)+g_N-\ind_{\{k=2\}}.
\]
In genus zero, any $d_k(N)$ distinct members form a basis, provided $d_k(N)\le2(k-1)$. The remaining rows follow from Proposition~\ref{proposition_real-periods} and Theorems~\ref{theorem_cyclic-periods} and~\ref{theorem_AL-periods}, using the dimension calculations in Appendix~\ref{appendix_dimension-checks}. In genus zero, Propositions~\ref{proposition_cyclic-selection} and~\ref{proposition_AL-selection} provide the stated basis selections. The weight-two cases use the corrections already included in the rank bounds.

For cyclic translation we use both orbits, taking $\ell=4$ for $N=16,32,48$, $\ell=6$ for $N=36$, $\ell=3$ for $N=18,27,45$, and $\ell=2$ otherwise. For the standard real locus and the Atkin--Lehner construction, we use the complete families on all $b_N$ arcs. The genus term is $g_N$ for the direct families and standard real arcs, $g_M$ with $M=N/\ell$ for cyclic translation, and $h_N$ for Atkin--Lehner translates.

The nonspanning entries follow from the dimension obstruction, either for the whole space or for an individual character component, or from the asymptotic rank formula \eqref{equation_real-period-asymptotic}. Here $k\gg1$ means all sufficiently large even weights at a fixed level.

\begin{table}[!htbp]
\caption{Scope of the complete families of critical periods for $1\le N<50$.}
\label{table_small-levels}
\centering
\small
\renewcommand{\arraystretch}{1.18}
\begin{tabular}{@{}p{0.19\textwidth}p{0.34\textwidth}p{0.41\textwidth}@{}}
\hline
Construction & Levels $N$ & Consequence \\ \hline
Periods at $0$
&
$1,2,3,4,5,6,7,8,9,11$
&
Spanning for every $k$.
\\
Periods at $0,1/2$
&
$10$
&
Genus zero; spanning and basis selection for every $k$.
\\
Periods at $0,1/2$
&
$14$
&
Genus one; spanning for every $k$.
\\
Periods at $0,1/2$
&
$22,26,30,34,38,42,46$
&
Does not span for $k\gg1$.
\\ \hline
Standard real arcs
&
$22,26,30,40$
&
Spanning for every $k$.
\\
Standard real arcs
&
$42,44$
&
$d_k(N)-1\le\rank\mathcal P\le d_k(N)$ for every $k$.
\\
Standard real arcs
&
$34,38,46$
&
Rank $4k+O_N(1)$; nonspanning for $k\gg1$.
\\ \hline
Cyclic translation
&
$8,12,16,18,20,24,$
$27,32,36,48$
&
Genus zero; spanning and basis selection for every $k$.
\\
Cyclic translation
&
$28,45$
&
Genus one; spanning for every $k$.
\\
Cyclic translation
&
$40,44$
&
Does not span for $k\gg1$.
\\ \hline
Atkin--Lehner
&
$5,7,9,11,13,15,$
$17,19,21,23,35$
&
Genus zero; spanning and basis selection for every $k$.
\\
Atkin--Lehner
&
$33$
&
Genus one; $d_k(N)-1\le\rank\mathcal P\le d_k(N)$ for $k\ge4$.
For $k=2$, the rank is $2<d_2(33)=3$.
\\
Atkin--Lehner
&
$25,27,29,31,37,39,$
$41,43,45,47,49$
&
Does not span for $k\gg1$.
\\ \hline
\end{tabular}
\end{table}

\appendix
\section{Dimension checks for Table~\ref{table_small-levels}}\label{appendix_dimension-checks}

Let $c_N,e_2(N),e_3(N)$ denote the numbers of cusps and elliptic points of orders two and three on $X_0(N)$. The standard dimension formula \cite[Proposition~6.1]{stein2007modular} is
\[
d_k(N)=(k-1)(g_N-1)+(k/2-1)c_N
       +e_2(N)\lfloor k/4\rfloor+e_3(N)\lfloor k/3\rfloor
\quad(k\ge4),
\]
with $d_2(N)=g_N$ in weight two. Combining these dimension formulas with the genus formula, we obtain
$d_k(N)+g_N-\ind_{\{k=2\}}\le k\mu_N/12-c_N+1$.
For $N=2,3,4,5,6,7,8,9,11$, we have $\mu_N\le12$ and $c_N\ge2$, so the threshold is at most $k-1$. At level one, $d_k(1)\le k/2-1$, and the $k/2-1$ critical indices with $n+1<k/2$ suffice by Theorem~\ref{theorem_basic}. For the two-path family at levels $10$ and $14$, the dimensions for $k\ge4$ are
\[
d_k(10)=k-3+2\lfloor k/4\rfloor,\qquad d_k(14)=2k-4.
\]
Their genera are zero and one, respectively, so $d_k(N)+g_N\le2(k-1)$. In weight two the dimensions are zero and one, and the corrected criterion also holds.

For the standard real arcs, we use the dimension formula and \eqref{equation_real-arc-count} to compare the spanning thresholds with the number of available periods:
\[
d_k(N)+g_N-b_N(k-1)=
\begin{cases}
1-k,&N=22,\\
1-k+2\lfloor k/4\rfloor,&N=26,\\
1-2k,&N=30,40,\\
1,&N=42,44
\end{cases}
\qquad(k\ge4).
\]
For the same levels in the same order, the corresponding weight-two differences $2g_N-1-b_N$ are $-1,-1,-3,-3,1,1$. By Proposition~\ref{proposition_real-periods}, the complete family therefore spans at $22,26,30,40$ and has rank at least $d_k(N)-1$ at $42,44$.

For cyclic translation, use the values of $\ell$ specified in Section~\ref{subsection_small-levels}, with $M=N/\ell$. A cusp $a/c$ has width $w=M/\gcd(c^2,M)$ and multiplier
\[
\chi_r(\gamma_{a/c})=\zeta^{r a^2w(1-acw)}.
\]
Let $b_r$ be the sum of the least positive cusp exponents, as in the verification of Example~\ref{remark_period-level27}. In genus zero, the valence calculation bounds the dimension by
\[
d_r\le \left[\frac{k\mu_M}{12}+1-b_r\right]_+.
\]
Here the nonnegative contributions from elliptic points have been omitted. Computing the cusp exponent sums from the multipliers above, we find the following minima, together with the group indices:
\[
(\mu_M,\min_r b_r)=
\begin{cases}
(6,2),&N=8,16,\\
(12,3),&N=12,18,\\
(12,2),&N=27,32,36,\\
(18,3),&N=20,\\
(24,4),&N=24,48.
\end{cases}
\]
Thus $d_r\le2(k-1)$ for every $r$ and every even $k\ge2$ in all these cases. Each level satisfies the spanning criterion in Theorem~\ref{theorem_cyclic-periods}, proving the genus-zero cyclic translation assertions in Table~\ref{table_small-levels}.

For $N=28,45$, the base curves have genus one, index $24$, and no elliptic points. Their cusp exponent sums are $(4,3)$ and $(4,3,3)$, respectively. For $k\ge4$, the cusp-form line bundles have positive degrees $2k-b_r$. On a genus-one curve, the Riemann--Roch formula identifies these degrees with the dimensions of the corresponding spaces of sections:
\[
\begin{aligned}
N=28:&\quad(d_0,d_1)=(2k-4,2k-3),\\
N=45:&\quad(d_0,d_1,d_2)=(2k-4,2k-3,2k-3).
\end{aligned}
\]
At $k=2$, the trivial-character space is the one-dimensional space of holomorphic differentials; every other component has a line bundle of degree one and hence dimension one. Thus the largest threshold $d_r+1-\ind_{\{k=2,r=0\}}$ equals $2(k-1)$ in every even weight. At $N=40,44$, however, the corresponding component has dimension $d_1=3k-4$ for $k\ge4$, which exceeds the available $2(k-1)$ moments.

For the Atkin--Lehner rows, we express the dimensions of the two sign components using the full Fricke trace formula \cite[Section~3.2, Proposition~3.1]{martin2023root}:
\[
\operatorname{tr}(W_N)=(-1)^{k/2}A_N+\ind_{\{k=2\}},\qquad
d_\pm=\frac{d_k(N)\pm(A_N-\ind_{\{k=2\}})}2,
\]
with the constants listed below. The last row records the smaller of $q_+-d_+$ and $q_--d_-$ at weight four, using all $b_N$ paths:
\[
\begin{array}{c|rrrrrrrrrrr}
N&5&7&9&11&13&15&17&19&21&23&35\\ \hline
A_N&1&1&1&2&1&2&2&2&2&3&4\\
\min(q_+-d_+,q_--d_-)\big|_{k=4}
&2&2&2&2&1&3&1&1&2&0&1
\end{array}
\]
At weight two, $h_N=d_-=(g_N-A_N+1)/2=0$ and $d_+=g_N=A_N-1\le b_N$. Together with the weight-four comparisons above, this verifies spanning in weights two and four. To pass to higher weights, applying the standard dimension formula at successive even weights, we obtain
\[
d_{k+2}(N)-d_k(N)
\le2g_N-2+c_N+e_2(N)+e_3(N)\le2b_N
\qquad(k\ge4).
\]
Hence each $d_\pm$ increases by at most $b_N$ when the weight increases by two, while each $q_\pm$ increases by exactly $b_N$. The spanning inequalities therefore hold in every even weight.

At $N=33$, the trace constant is $A_{33}=2$ and $h_{33}=1$. For $k\ge4$,
\[
(d_+,d_-)=(2k-2,2k-4),\qquad(q_+,q_-)=(2k,2k-4),
\]
so $d_k(33)-1\le\rank\mathcal P\le d_k(33)$ by \eqref{equation_AL-period-rank}. At weight two, the dimensions $(2,1)$ and moment counts $(4,0)$ make both rank bounds equal to two.

Finally, the two-path families at levels $22,26,30,34,38,42,46$ cannot span in sufficiently large weights: here $\mu_N>24$, so $d_k(N)\sim k\mu_N/12>2(k-1)$ for sufficiently large $k$. For the standard real arcs at $34,38,46$, the values are $b_N=4$ and $\mu_N=54,60,72$, respectively. Their rank is therefore $4k+O_N(1)$ by \eqref{equation_real-period-asymptotic}. At every level in the Atkin--Lehner nonspanning row, $\mu_N>12b_N$, so the same formula places the rank strictly below the dimension for sufficiently large weights.

\providecommand{\bysame}{\leavevmode\hbox to3em{\hrulefill}\thinspace}
\providecommand{\MR}{\relax\ifhmode\unskip\space\fi MR }
\providecommand{\MRhref}[2]{%
  \href{http://www.ams.org/mathscinet-getitem?mr=#1}{#2}
}
\providecommand{\href}[2]{#2}

\end{document}